\documentclass[journal]{IEEEtran}%
\usepackage{amssymb}
\usepackage{amsmath}
\usepackage{amsfonts}
\usepackage{graphicx}%
\providecommand{\U}[1]{\protect\rule{.1in}{.1in}}
\newtheorem{theorem}{Theorem}

\newtheorem{lemma}[theorem]{Lemma}

\begin{document}

\title{In-homogeneous Virus Spread in Networks}
\author{Piet Van Mieghem and~Jasmina Omic\thanks{Delft University of Technology,
Faculty of Electrical Engineering, Mathematics and Computer Science, P.O Box
5031, 2600 GA Delft, The Netherlands. Email: \{P.VanMieghem,
J.S.Omic\}@ewi.tudelft.nl. This work of 2008 was a Delft University of
Technology, report2008081 on
http://www.nas.ewi.tudelft.nl/people/Piet/TUDelftReports.html}}
\maketitle

\begin{abstract}
Our $N$-intertwined mean-field approximation (NIMFA)
\cite{PVM_ToN_VirusSpread} for virus spread in any network with $N$ nodes is
extended to a full heterogeneous setting. The metastable steady-state nodal
infection probabilities are specified in terms of a generalized Laplacian,
that possesses analogous properties as the classical Laplacian in graph
theory. The critical threshold that separates global network infection from
global network health is characterized via an $N$ dimensional vector that
makes the largest eigenvalue of a modified adjacency matrix equal to unity.
Finally, the steady-state infection probability of node $i$ is convex in the
own curing rate $\delta_{i}$, but can be concave in the curing rates
$\delta_{j}$ of the other nodes $1\leq j\neq i\leq N$ in the network.

\end{abstract}

\begin{keywords}
Virus spread, epidemic threshold, generalized Laplacian
\end{keywords}

\section{Introduction}

This paper generalizes our $N$-Intertwined Mean-Field Approximation (NIMFA)
for virus spread in networks, presented in \cite{PVM_ToN_VirusSpread} and
\cite[Chapter 17]{PVM_PAComplexNetsCUP}, to a heterogeneous setting.
\emph{Heterogeneity} rather than \emph{homogeneity} abounds in real networks.
For example, in data communications networks, the transmission capacity, age,
performance, installed software, security level and other properties of
networked computers are generally different. Social and biological networks
are very diverse: a population often consists of a mix of weak and strong, or
old and young species or of completely different types of species. The network
topology for transport by airplane, car, train, ship is different. Many more
examples can be added illustrating that homogeneous networks are the exception
rather than the rule. This diversity in the \textquotedblleft
nodes\textquotedblright\ and \textquotedblleft links\textquotedblright\ of
real networks will thus likely affect the spreading pattern of viruses, that
are here understood as malicious challenges of a network.

NIMFA approximates the continuous-time Markov susceptible-infected-susceptible
(SIS) epidemic process on a network with $N$ nodes, that was earlier
considered by Ganesh \emph{et al. }\cite{Ganesh_2005} and by Wang \emph{et
al.} \cite{YWang2003} in discrete-time. Each node in the network is either
infected or healthy. In a heterogeneous setting, an infected node $i$ can
infect its neighbors with an infection rate $\beta_{i}$, but it is cured with
curing rate $\delta_{i}$. Once cured and healthy, the node is again prone to
the virus. Both infection and curing processes are independent.

Previously in \cite{PVM_ToN_VirusSpread}, only a homogeneous virus spread was
investigated, where all infection rates $\beta_{i}=\beta$ and all curing rates
$\delta_{i}=\delta$ were the same for each node. We believe that the extension
to a full heterogeneous setting is, perhaps, the best SIS model that we can
achieve. The exact Markovian model, described and analyzed in
\cite{PVM_ToN_VirusSpread}, has $2^{N}$ states, which makes it infeasible to
compute for realistic sizes of networks. Moreover, the exact Markovian model
possesses as steady-state the overall healthy state, which is an absorbing
state, that is, unfortunately, only reached after an extreme and
unrealistically long time. The heterogeneous NIMFA makes one approximation, a
mean field approximation as shown in Section \ref{sec_N_intertwinedMC} and in
\cite{PVM_ToN_VirusSpread}, that results in a set of $N$ non-linear equations.
Hence, NIMFA trades computational feasibility, a reduction of $2^{N}$ linear
equations to $N$ non-linear ones, at the expense of exactness. The last point,
the accuracy of NIMFA is shown in \cite{PVM_ToN_VirusSpread} (and further in
\cite{Jasmina_variance2008}) to be overall remarkably good, with a worst case
performance near the critical threshold, which is a realistic and observable
artifact of the metastable steady-state that does not exist in the exact
Markovian steady-state. Below the critical epidemic threshold, infection
vanishes exponentially fast in time and above the critical threshold the
network stays infected to a degree determined by the effective infection
vector $\tau$, with components $\tau_{i}=\frac{\beta_{i}}{\delta_{i}}$.

A major new insight is that the metastable steady-state can be written in
terms of a generalized Laplacian matrix that bears similar deep properties as
the Laplacian matrix of a graph (see e.g. \cite{Biggs}, \cite{Cvetkovic} and
\cite{PVM_graphspectra}). In a heterogeneous setting, the critical threshold
is characterized by an effective infection vector, instead of one scalar in
the homogeneous case equal to $\tau_{\text{hom;}c}=\frac{1}{\lambda_{\max
}\left(  A\right)  }$, where $\lambda_{\max}\left(  A\right)  $ is the largest
eigenvalue of the adjacency matrix $A$ of the graph. This critical vector
determines a critical surface in the $N$-dimensional space spanned by the
vector components $\tau_{1},\ldots,\tau_{N}$. We also prove that the
steady-state infection probability $v_{i\infty}$ of node $i$ is convex in the
curing rate $\delta_{i}$, given all other curing rates $\delta_{j}$ are the
same.\textbf{ }This convexity result is applied in a virus protection game
played by the individual and selfish nodes in a network
\cite{PVM_Jasmina_Game_protection_INFOCOM2009}.

\section{$N$-intertwined continuous Markov chains with $2$ states}

\label{sec_N_intertwinedMC}This section extends the homogeneous NIMFA in
\cite{PVM_ToN_VirusSpread} to a heterogeneous setting. Although analogous to
the corresponding section in \cite{PVM_ToN_VirusSpread}, its inclusion makes
this paper self-contained.

By separately observing each node, we will model the virus spread in a
bi-directional network specified by a symmetric adjacency matrix $A$. Every
node $i$ at time $t$ in the network has two states: infected with probability
$\Pr[X_{i}\left(  t\right)  =1]$ and healthy with probability $\Pr
[X_{i}\left(  t\right)  =0]$. At each moment $t$, a node can only be in one of
two states, thus $\Pr[X_{i}(t)=1]+\Pr[X_{i}(t)=0]=1$. If we apply Markov
theory, the infinitesimal generator $Q_{i}\left(  t\right)  $ of this
two-state continuous Markov chain is,%
\[
Q_{i}\left(  t\right)  =\left[
\begin{array}
[c]{cc}%
-q_{1;i} & q_{1;i}\\
q_{2;i} & -q_{2;i}%
\end{array}
\right]
\]
with $q_{2;i}=\delta_{i}$ and
\[
q_{1;i}=\sum_{j=1}^{N}\beta_{j}a_{ij}1_{\left\{  X_{j}(t)=1\right\}  }%
\]
where the indicator function $1_{x}=1$ if the event $x$ is true else it is
zero. The coupling of node $i$ to the rest of the network is described by an
infection rate $q_{1;i}$ that is a random variable, which essentially makes
the process doubly stochastic. This observation is crucial. For, using the
definition of the infinitesimal generator \cite[p. 181]%
{PVM_PerformanceAnalysisCUP},%
\[
\Pr[\left.  X_{i}(t+\Delta t)=1\right\vert X_{i}\left(  t\right)
=0]=q_{1;i}\Delta t+o(\Delta t)
\]
the continuity and differentiability shows that this process is not Markovian
anymore. The random nature of $q_{1;i}$ is removed by an additional
conditioning to all possible combinations of rates, which is equivalent to
conditioning to all possible combinations of the states $X_{j}(t)=1$ (and
their complements $X_{j}(t)=0$) of the neighbors of node $i$. Hence, the
number of basic states dramatically increases. Eventually, after conditioning
each node in such a way, we end up with a $2^{N}$-- state Markov chain,
studied in \cite{PVM_ToN_VirusSpread}.

Instead of conditioning, we replace the actual, random infection rate by an
effective or average infection rate, which is basically a mean field
approximation,
\begin{equation}
E\left[  q_{1;i}\right]  =E\left[  \sum_{j=1}^{N}\beta_{j}a_{ij}1_{\left\{
X_{j}(t)=1\right\}  }\right]  \label{mean_field_approximation}%
\end{equation}
In general, we may take the expectation over the rates $\beta_{i}$, the
network topology via the matrix $A$ and the states $X_{j}(t)$. Since we assume
that both the infection rates $\beta_{i}$ and the network are constant and
given, we only average over the states. Using $E\left[  1_{x}\right]
=\Pr\left[  x\right]  $ (see e.g. \cite{PVM_PerformanceAnalysisCUP}), we
replace $q_{1;i}$ by
\[
E\left[  q_{1;i}\right]  =\sum_{j=1}^{N}\beta_{j}a_{ij}\Pr[X_{j}(t)=1]
\]
which results in an effective infinitesimal generator,%
\[
\overline{Q_{i}(t)}=\left[
\begin{array}
[c]{cc}%
-E\left[  q_{1;i}\right]  & E\left[  q_{1;i}\right] \\
\delta_{i} & -\delta_{i}%
\end{array}
\right]
\]

The effective $\overline{Q_{i}(t)}$ allows us to proceed with Markov theory.
Denoting $v_{i}\left(  t\right)  =\Pr[X_{i}(t)=1]$ and recalling that
$\Pr[X_{i}(t)=0]=1-v_{i}\left(  t\right)  $, the Markov differential equation
\cite[(10.11) on p. 208]{PVM_PAComplexNetsCUP} for state $X_{i}(t)=1$ turns
out to be non-linear%
\begin{equation}
\frac{dv_{i}\left(  t\right)  }{dt}=\sum_{j=1}^{N}\beta_{j}a_{ij}v_{j}\left(
t\right)  -v_{i}\left(  t\right)  \left(  \sum_{j=1}^{N}\beta_{j}a_{ij}%
v_{j}\left(  t\right)  +\delta_{i}\right)  \label{Non_linear_dvg_state_Xi=1}%
\end{equation}
Each node obeys a differential equation as (\ref{Non_linear_dvg_state_Xi=1}),%
\[
\hspace{-0.3cm}\left\{
\begin{array}
[c]{c}%
\frac{dv_{1}\left(  t\right)  }{dt}=\sum_{j=1}^{N}\beta_{j}a_{1j}%
v_{j}\!\left(  t\right)  \!-v_{1}\!\left(  t\right)  \!\left(  \sum_{j=1}%
^{N}\beta_{j}a_{1j}v_{j}\!\left(  t\right)  +\delta_{1}\right) \\
\frac{dv_{2}\left(  t\right)  }{dt}=\sum_{j=1}^{N}\beta_{j}a_{2j}%
v_{j}\!\left(  t\right)  \!-v_{2}\!\left(  t\right)  \!\left(  \sum_{j=1}%
^{N}\beta_{j}a_{2j}v_{j}\!\left(  t\right)  +\delta_{2}\right) \\
\vdots\\
\frac{dv_{N}\left(  t\right)  }{dt}=\sum_{j=1}^{N}\beta_{j}a_{Nj}%
v_{j}\!\left(  t\right)  \!-v_{N}\!\left(  t\right)  \!\left(  \sum_{j=1}%
^{N}\beta_{j}a_{Nj}v_{j}\!\left(  t\right)  +\delta_{N}\right)
\end{array}
\right.
\]
Written in matrix form, with
\[
V\left(  t\right)  =\left[
\begin{array}
[c]{cccc}%
v_{1}\left(  t\right)  & v_{2}\left(  t\right)  & \cdots & v_{N}\left(
t\right)
\end{array}
\right]  ^{T}%
\]
we arrive at%
\begin{equation}
\frac{dV\left(  t\right)  }{dt}=A\text{diag}\left(  \beta_{j}\right)  V\left(
t\right)  -\text{diag}\left(  v_{i}\left(  t\right)  \right)  \left(
A\text{diag}\left(  \beta_{j}\right)  V\left(  t\right)  +C\right)
\label{non_linear_matrix_evolution_diffequation}%
\end{equation}
where diag$\left(  v_{i}\left(  t\right)  \right)  $ is the diagonal matrix
with elements $v_{1}\left(  t\right)  ,v_{2}\left(  t\right)  ,\ldots
\,,v_{N}\left(  t\right)  $ and the curing rate vector is $C=\left(
\delta_{1},\delta_{2},\ldots,\delta_{N}\right)  $.

We note that $A$diag$\left(  \beta_{i}\right)  $ is, in general and opposed to
the homogeneous setting, not symmetric anymore, unless $A$ and diag$\left(
\beta_{i}\right)  $ commute, in which case the eigenvalue $\lambda_{i}\left(
A\text{diag}\left(  \beta_{i}\right)  \right)  =\lambda_{i}\left(  A\right)
\beta_{i}$ and both $\beta_{i}$ and $\lambda_{i}\left(  A\right)  $ have a
same eigenvector $x_{i}$.

\section{General in-homogenous steady-state}

\label{sec_in_homogeneous_curing}

\subsection{The steady-state equation}

The metastable steady-state follows from
(\ref{non_linear_matrix_evolution_diffequation}) as%
\[
A\text{diag}\left(  \beta_{i}\right)  V_{\infty}-\text{diag}\left(
v_{i\infty}\right)  \left(  A\text{diag}\left(  \beta_{i}\right)  V_{\infty
}+C\right)  =0
\]
where $V_{\infty}=\lim_{t\rightarrow\infty}V\left(  t\right)  $. We define the
vector
\begin{equation}
w=A\text{diag}\left(  \beta_{i}\right)  V_{\infty}+C \label{def_w}%
\end{equation}
and write the stead-state equation as%
\[
w-C=\text{diag}\left(  v_{i\infty}\right)  w
\]
or%
\[
\left(  I-\text{diag}\left(  v_{i\infty}\right)  \right)  w=C
\]
Ignoring extreme virus spread conditions (the absence of curing ($\delta
_{i}=0$) and an infinitely strong infection rate $\beta_{i}\rightarrow\infty
$), then the infection probabilities $v_{i\infty}$ cannot be one such that the
matrix $\left(  I-\text{diag}\left(  v_{i\infty}\right)  \right)  =$
diag$\left(  1-v_{i\infty}\right)  $ is invertible. Hence,%
\[
w=\text{diag}\left(  \frac{1}{1-v_{i\infty}}\right)  C
\]
Invoking the definition (\ref{def_w}) of $w$, we obtain%
\begin{align}
A\text{diag}\left(  \beta_{i}\right)  V_{\infty}  &  =\text{diag}\left(
\frac{v_{i\infty}}{1-v_{i\infty}}\right)  C\nonumber\\
&  =\text{diag}\left(  \frac{\delta_{i}}{1-v_{i\infty}}\right)  V_{\infty}
\label{steady_state_equation}%
\end{align}
The $i$-th row of (\ref{steady_state_equation}) yields the nodal steady state
equation,%
\begin{equation}
\sum_{j=1}^{N}a_{ij}\beta_{j}v_{j\infty}=\frac{v_{i\infty}\delta_{i}%
}{1-v_{i\infty}} \label{steady_state_node_i_equation}%
\end{equation}
Let $\widetilde{V}_{\infty}=$ diag$\left(  \beta_{i}\right)  V_{\infty}$ and
the effective spreading rate for node $i$, $\tau_{i}=\frac{\beta_{i}}%
{\delta_{i}}$, then we arrive at%
\begin{equation}
\mathcal{Q}\left(  \frac{1}{\tau_{i}\left(  1-v_{i\infty}\right)  }\right)
\widetilde{V}_{\infty}=0 \label{eigenvector_state_state_solution}%
\end{equation}
where the symmetric matrix
\begin{align}
\mathcal{Q}\left(  q_{i}\right)   &  \mathcal{=}\text{ diag}\left(
q_{i}\right)  -A\label{def_generalized_Laplacian}\\
&  =\text{diag}\left(  q_{i}-d_{i}\right)  +Q\nonumber
\end{align}
can be interpreted as a generalized Laplacian\footnote{All eigenvalues of the
Laplacian $Q=\Delta-A$ in a connected graph are positive, except for the
smallest one that is zero. Hence, $Q$ is positive semi-definite. Much more
properties of the Laplacian $Q$ are found e.g. in \cite{Biggs} and
\cite{Cvetkovic}.}, because $\mathcal{Q}\left(  d_{i}\right)  =Q=\Delta-A$,
where $\Delta=$ diag$\left(  d_{i}\right)  $. The observation that the
non-linear set of steady-state equations can be written in terms of the
generalized Laplacian $\mathcal{Q}\left(  q_{i}\right)  $ is fortunate,
because, as will be shown in Section \ref{sec_generalized_Laplacian}, the
powerful theory of the \textquotedblleft normal\textquotedblright\ Laplacian
$Q$ applies.

The modified steady-state vector $\widetilde{V}_{\infty}$ is orthogonal to
each row (or, by symmetry, each column) vector of $\mathcal{Q}\left(  \frac
{1}{\tau_{i}\left(  1-v_{i\infty}\right)  }\right)  $. A non-zero modified
steady-state vector $\widetilde{V}_{\infty}$ is thus only possible provided
$\det\mathcal{Q}\left(  \frac{1}{\tau_{i}\left(  1-v_{i\infty}\right)
}\right)  =0$. In other words, the generalized Laplacian $\mathcal{Q}\left(
\frac{1}{\tau_{i}\left(  1-v_{i\infty}\right)  }\right)  $ should have a zero
eigenvalue with the modified steady-state vector $\widetilde{V}_{\infty}$ as
corresponding eigenvector. Since the vectors $B=\left(  \beta_{1},\beta
_{2},\ldots,\beta_{N}\right)  $ and $C=\left(  \delta_{1},\delta_{2}%
,\ldots,\delta_{N}\right)  $ are given, the non-linear eigenvector problem
(\ref{eigenvector_state_state_solution}) has, in general, a solution that
cannot simply be recast to the homogeneous case where $B=\beta u$ and
$C=\delta u$ (or $\beta_{i}=\beta$ and $\delta_{i}=\delta$ for all $1\leq
i\leq N$) in which the all-one vector $u=\left(  1,1,\ldots,1\right)  $.

\subsection{The generalized Laplacian $\mathcal{Q}\left(  q_{i}\right)  $}

\label{sec_generalized_Laplacian}Since $\mathcal{Q}\left(  q_{i}\right)  $ is
symmetric, all eigenvectors are orthogonal such that, with $\widetilde
{V}_{\infty}=$ diag$\left(  \beta_{i}\right)  V_{\infty}$,%
\begin{equation}
\sum_{j=1}^{N}\beta_{j}v_{j\infty}y_{j}=0
\label{general_eigenvector_orthogonality}%
\end{equation}
where $y$ is the eigenvector belonging to eigenvalue $\lambda\left(
\mathcal{Q}\left(  q_{i}\right)  \right)  \neq0$.

\begin{theorem}
\label{theorem_all_eigenvalues_generalized_Q_are_positive_except_1}If the
network $G$ is connected, all eigenvalues of $\mathcal{Q}\left(  q_{i}\right)
$ are positive, except for the smallest one $\lambda_{N}\left(  \mathcal{Q}%
\right)  =0$.
\end{theorem}

\textbf{Proof:} The theorem is a consequence of the Perron-Frobenius Theorem
(see e.g. \cite{GantmacherII}) for a non-negative, irreducible matrix. Indeed,
consider the non-negative matrix $q_{\max}I-\mathcal{Q}\left(  q_{i}\right)
$, where $q_{\max}=\max_{1\leq i\leq N}q_{i}$, whose eigenvalues are $\xi
_{k}=q_{\max}-\lambda_{k}\left(  \mathcal{Q}\right)  $ for $1\leq k\leq N$.
Since $G$ is connected, then $q_{\max}I-\mathcal{Q}\left(  q_{i}\right)  $ is
irreducible and the Perron-Frobenius Theorem states that the largest
eigenvalue $r=\max_{1\leq k\leq N}\xi_{k}$ of $q_{\max}I-\mathcal{Q}\left(
q_{i}\right)  $ is positive and simple and the corresponding eigenvector
$x_{r}$ has positive components. Hence, $\mathcal{Q}\left(  q_{i}\right)
x_{r}=\left(  q_{\max}-r\right)  x_{r}$. Since eigenvectors of a symmetric
matrix are orthogonal while $\widetilde{V}_{\infty}^{T}x_{r}>0$, $x_{r}$ must
be proportional to $\widetilde{V}_{\infty}$, and thus $q_{\max}=r$. Since
there is only one such eigenvector $x_{r}$ and since the eigenvalue
$r>q_{\max}-\lambda_{k}\left(  \mathcal{Q}\right)  $ for all $k$ (except that
$k$ for which $\lambda_{k}\left(  \mathcal{Q}\right)  =0$, which is thus the
smallest eigenvalue), all other eigenvalues of $\mathcal{Q}\left(
q_{i}\right)  $ must exceed zero.\hfill$\square\medskip$

If the graph $G$ is disconnected which means that $A$ is reducible
\cite{PVM_PerformanceAnalysisCUP}, the Theorem
\ref{theorem_all_eigenvalues_generalized_Q_are_positive_except_1} still
applies (see e.g. \cite{GantmacherII}), however, under the slightly weakened
form that $x_{r}$ has non-negative components (instead of positive, hence,
zero components can occur) and that the largest eigenvalue $r$ is non-zero
(not necessarily strict positive). The consequence is that more than one zero
eigenvalue can occur. From the point of virus spread, we may ignore
disconnected graphs, because the theory can be applied to each connected
component (cluster) of the network $G$. The symmetry of $\mathcal{Q}\left(
q_{i}\right)  $ implies that all eigenvalues are real and can be ordered. By
Theorem \ref{theorem_all_eigenvalues_generalized_Q_are_positive_except_1}, we
have
\[
0=\lambda_{N}\left(  \mathcal{Q}\right)  \leq\lambda_{N-1}\left(
\mathcal{Q}\right)  \leq\ldots\leq\lambda_{1}\left(  \mathcal{Q}\right)
\]
Gerschgorin's theorem \cite[p. 71-75]{Wilkinson} indicates that the
eigenvalues of $\mathcal{Q}\left(  q_{i}\right)  $ are centered around $q_{i}$
with radius equal to the degree $d_{i}$, i.e. an eigenvalue $\lambda$ of
$\mathcal{Q}\left(  q_{i}\right)  $ lies in an interval $\left\vert
\lambda-q_{k}\right\vert \leq d_{k}$ for some $1\leq k\leq N$. Thus, there is
an eigenvalue $\lambda$ of $\mathcal{Q}\left(  q_{i}\right)  $ that obeys
\[
q_{k}-d_{k}\leq\lambda\leq d_{k}+q_{k}%
\]
A solution of (\ref{eigenvector_state_state_solution}) requires that at least
one eigenvalue of $\mathcal{Q}\left(  q_{i}\right)  $ is zero, while Theorem
\ref{theorem_all_eigenvalues_generalized_Q_are_positive_except_1} states that
there is only one zero eigenvalue. Hence, precisely one, say the $j$-th, of
the Gerschgorin line segments that contain the eigenvalue $\lambda_{N}\left(
\mathcal{Q}\right)  =0$, must obey $q_{j}\leq d_{j}$ to have a non-zero
solution of (\ref{eigenvector_state_state_solution}). However, more
Gerschgorin segments may obey $q_{k}-d_{k}\leq0$. This couples $\frac{1}%
{\tau_{j}\left(  1-v_{j\infty}\right)  }\leq d_{j}$ for at least one $j$
component and shows that, when $v_{j\infty}\rightarrow1$, there must hold that
$\tau_{j}\rightarrow\infty$. Hence, for at least one component $j$, there
holds that%
\[
0<v_{j\infty}\leq1-\frac{1}{\tau_{j}d_{j}}%
\]
where the lower bound follows, by the Perron-Frobenius Theorem, from the fact
that the network $G$ is connected. This shows that there is a critical bound
on $\tau_{j}>\frac{1}{d_{j}}$ for at least one component of $\tau$. The
critical threshold on the $\tau$-vector is further explored in Section
\ref{sec_critical_threshold}, while Section \ref{sec_critical_vector_KN}
applies the theory to the complete graph.

We also know that trace$\left(  \mathcal{Q}\left(  q_{i}\right)  \right)
=\sum_{k=1}^{N}\lambda_{k}\left(  \mathcal{Q}\right)  $. Thus, with
$\lambda_{N}\left(  \mathcal{Q}\right)  =0$,%
\[
\sum_{k=1}^{N-1}\lambda_{k}\left(  \mathcal{Q}\right)  =\sum_{i=1}^{N}\frac
{1}{\tau_{i}\left(  1-v_{i\infty}\right)  }%
\]
In addition, since
\begin{align*}
\text{trace}\left(  \mathcal{Q}^{2}\left(  q_{i}\right)  \right)   &
=\text{trace}\left(  \text{diag}\left(  q_{i}^{2}\right)  \right)
+\text{trace}\left(  A^{2}\right) \\
&  =\sum_{i=1}^{N}\frac{1}{\tau_{i}^{2}\left(  1-v_{i\infty}\right)  ^{2}}+2L
\end{align*}
we have that%
\[
\sum_{k=1}^{N-1}\lambda_{k}^{2}\left(  \mathcal{Q}\right)  =\sum_{i=1}%
^{N}\frac{1}{\tau_{i}^{2}\left(  1-v_{i\infty}\right)  ^{2}}+2L
\]

Right multiplication of (\ref{steady_state_equation}) by the all one-vector
$u^{T}=(1,1,\ldots,1)$ yields%
\[
u^{T}A\text{diag}\left(  \beta_{i}\right)  V_{\infty}=u^{T}\text{diag}\left(
\frac{\delta_{i}}{1-v_{i\infty}}\right)  V_{\infty}%
\]
With $u^{T}A=D^{T}=\left(  d_{1},d_{2},\ldots,d_{N}\right)  $, the degree
vector, we have%
\[
\left(  u^{T}\text{diag}\left(  \frac{\delta_{i}}{1-v_{i\infty}}\right)
-D^{T}\text{diag}\left(  \beta_{i}\right)  \right)  V_{\infty}=0
\]
or\footnote{The result (\ref{steady_state_constraint_vgl}) also follows by
adding all rows in (\ref{eigenvector_state_state_solution})
\[
\mathcal{Q}\left(  q_{i}\right)  \widetilde{V}_{\infty}=\text{diag}\left(
q_{i}-d_{i}\right)  \widetilde{V}_{\infty}+Q\widetilde{V}_{\infty}%
\]
and using the basic fact that the row sum of the Laplacian $Q$ is zero.}
\begin{equation}
\sum_{j=1}^{N}\left(  \frac{1}{\tau_{j}\left(  1-v_{j\infty}\right)  }%
-d_{j}\right)  \beta_{j}v_{j\infty}=0 \label{steady_state_constraint_vgl}%
\end{equation}
Similarly as deduced from Gershgorin's theorem, this sum shows that, at least
one $j$ term should be negative (because $\beta_{j}v_{j\infty}\geq0$), i.e.
$d_{j}\geq\frac{1}{\tau_{j}\left(  1-v_{j\infty}\right)  }$. Also, in view of
(\ref{general_eigenvector_orthogonality}), the vector $y$ with components
$y_{j}=\frac{1}{\tau_{j}\left(  1-v_{j\infty}\right)  }-d_{j}$ is a linear
combination of eigenvectors of $\mathcal{Q}\left(  \frac{1}{\tau_{i}\left(
1-v_{i\infty}\right)  }\right)  $ belonging to a non-zero eigenvalue. In
general, however, the vector $y$ is not an eigenvector of $\mathcal{Q}\left(
\frac{1}{\tau_{i}\left(  1-v_{i\infty}\right)  }\right)  $.

\begin{lemma}
\label{lemma_generalized_laplacian_q*>q_positive_def}If $q_{i}^{\ast}>q_{i}$
for all $1\leq i\leq N$, then $\mathcal{Q}\left(  q_{i}^{\ast}\right)  $ is
positive definite.
\end{lemma}

\textbf{Proof:} For any non-zero vector $x$, consider the quadratic form%
\[
x^{T}\mathcal{Q}\left(  q_{i}^{\ast}\right)  x=x^{T}\mathcal{Q}\left(
q_{i}\right)  x+x^{T}\text{diag}\left(  q_{i}^{\ast}-q_{i}\right)  x
\]
Theorem \ref{theorem_all_eigenvalues_generalized_Q_are_positive_except_1}
implies that $x^{T}\mathcal{Q}\left(  q_{i}\right)  x\geq0$, i.e. that
$\mathcal{Q}\left(  q_{i}\right)  $ is semi-definite. Since $q_{i}^{\ast
}>q_{i}$ for all $1\leq i\leq N$, $x^{T}$diag$\left(  q_{i}^{\ast}%
-q_{i}\right)  x>0$, which demonstrates the lemma.\hfill$\square\medskip$

Lemma \ref{lemma_generalized_laplacian_q*>q_positive_def} indicates that the
matrix $\mathcal{Q}\left(  \frac{1}{\tau_{i}\left(  1-v_{i\infty}\right)
^{2}}\right)  $, that appears in the definition (\ref{def_S}) of the matrix
$S$ in Section \ref{sec_convexity}\textbf{,} is positive definite, because
$\mathcal{Q}\left(  \frac{1}{\tau_{i}\left(  1-v_{i\infty}\right)  }\right)  $
defines the vector $V_{\infty}=\left(  v_{1\infty},v_{2\infty},\ldots
,v_{N\infty}\right)  $ via (\ref{eigenvector_state_state_solution}).

\subsection{The critical threshold}

\label{sec_critical_threshold}We known that the exact steady-state is
$V_{\infty}=0$, but the metastable steady-state (see
\cite{PVM_ToN_VirusSpread} for a deeper discussion) is characterized by a
second solution, the eigenvector of (\ref{eigenvector_state_state_solution}).

\begin{theorem}
\label{theorem_critical_threshold}The critical threshold is determined by
vectors $\tau_{c}=$ $\left(  \tau_{1c},\tau_{2c},\ldots,\tau_{Nc}\right)  $
that obey $\lambda_{\max}\left(  R\right)  =1$, where $\lambda_{\max}\left(
R\right)  $ is the largest eigenvalue of the symmetric matrix
\begin{equation}
R=\text{diag}\left(  \sqrt{\tau_{i}}\right)  A\text{diag}\left(  \sqrt
{\tau_{i}}\right)  \label{def_R}%
\end{equation}
whose corresponding eigenvector has positive components if the graph $G$ is connected.
\end{theorem}

\textbf{Proof:} At the critical threshold, the second, non-zero solution is
$V_{\infty}=\varepsilon x$, where $x$ is a vector with non-negative components
and where $\varepsilon$ is arbitrary small. This property allows us to
approximate the generalized Laplacian $\mathcal{Q}\left(  q\right)  $ as%
\begin{align*}
\mathcal{Q}\left(  \frac{1}{\tau_{i}\left(  1-v_{i\infty}\right)  }\right)
&  \mathcal{=}\text{ diag}\left(  \frac{\delta_{i}}{\beta_{i}\left(
1-\varepsilon x_{i}\right)  }\right)  -A\\
&  =\text{diag}\left(  \frac{\delta_{i}}{\beta_{i}}\right)  \left(
I-\varepsilon\text{diag}\left(  x_{i}\right)  \right)  -A+O\left(
\varepsilon^{2}\right)
\end{align*}
such that (\ref{eigenvector_state_state_solution}) becomes to first order in
$\varepsilon$%
\[
\mathcal{Q}\left(  \frac{1}{\tau_{i}}\right)  \text{diag}\left(  \beta
_{i}\right)  x=0
\]
which can be rewritten as an eigenvalue equation for the adjacency matrix,%
\[
\text{diag}\left(  \frac{1}{\delta_{i}}\right)  A\text{diag}\left(  \beta
_{i}\right)  x=x
\]
Hence, $x$ is the eigenvector of $\widetilde{A}=$ diag$\left(  \frac{1}%
{\delta_{i}}\right)  A$diag$\left(  \beta_{i}\right)  $ belonging to the
eigenvalue 1. Since $\widetilde{A}$ is a non-negative, irreducible matrix, the
Perron-Frobenius Theorem \cite[p. 451]{PVM_PerformanceAnalysisCUP} states that
$\widetilde{A}$ has a positive largest eigenvalue $\lambda_{\max}\left(
\widetilde{A}\right)  $ with a corresponding eigenvector whose elements are
all positive and that there is only one eigenvector of $\widetilde{A}$ with
non-negative components. Since any scaled vector $V_{\infty}=\varepsilon x$
must have non-negative components (because they represent scaled
probabilities), we find that $\lambda_{\max}\left(  \widetilde{A}\right)  =1$.
Hence, for the given vectors $B=\left(  \beta_{1},\beta_{2},\ldots,\beta
_{N}\right)  $ and $C=\left(  \delta_{1},\delta_{2},\ldots,\delta_{N}\right)
$, there are three possibilities:%
\[
\left\{
\begin{array}
[c]{cc}%
\lambda_{\max}\left(  \widetilde{A}\right)  <1 & \text{not infected network}\\
\lambda_{\max}\left(  \widetilde{A}\right)  =1 & \text{critical threshold}\\
\lambda_{\max}\left(  \widetilde{A}\right)  >1 & \text{infected network}%
\end{array}
\right.
\]
where the inequalities sign are deduced by relating the largest eigenvalue to
the norm of the matrix $\widetilde{A}$: higher eigenvalues correspond to a
larger norm (see e.g. \cite[Section A.3.1]{PVM_PerformanceAnalysisCUP}). Of
course, only in case $\lambda_{\max}\left(  \widetilde{A}\right)  =1$, the
eigenvector equation has a non-zero solution. If $\lambda_{\max}\left(
\widetilde{A}\right)  >1$, then the first order expansion is inadequate and
the full non-linear equation (\ref{eigenvector_state_state_solution}) needs to
be solved.

The first order expansion process has caused $\widetilde{A}$ to be not
symmetric, while $\mathcal{Q}\left(  \frac{1}{\tau_{i}\left(  1-v_{i\infty
}\right)  }\right)  $ is symmetric in general. Fortunately, there exist a
similarity transform $H=$ diag$\left(  \sqrt{\delta_{i}\beta_{i}}\right)  $
which symmetrizes $\widetilde{A}$,%
\[
R=H\widetilde{A}H^{-1}=\text{diag}\left(  \sqrt{\frac{\beta_{i}}{\delta_{i}}%
}\right)  A\text{diag}\left(  \sqrt{\frac{\beta_{i}}{\delta_{i}}}\right)
\]
and $R=R^{T}$ has the same real eigenvalues as $\widetilde{A}$ (see \cite[p.
438]{PVM_PerformanceAnalysisCUP}). The matrix $R$ also demonstrates that only
an effective rate per node, $\tau_{i}=\frac{\beta_{i}}{\delta_{i}}$, is
needed. Thus, the equation that characterizes the critical threshold is
\[
Ry=y
\]
where $y=Hx$. The eigenvalue $\lambda_{\max}\left(  \widetilde{A}\right)
=\lambda_{\max}\left(  R\right)  =1$ determines the critical vectors $\tau
_{c}=$ $\left(  \tau_{1c},\tau_{2c},\ldots,\tau_{Nc}\right)  $. In general,
there can be more than one critical vector because $\lambda_{\max}\left(
R\right)  =1$ is a map of $\mathbb{R}^{N}\rightarrow\mathbb{R}$.
\hfill$\square$

We remark that, since trace$\left(  R\right)  =$ trace$\left(  A\right)  =0$,
that $\lambda_{\max}\left(  R\right)  =\lambda_{1}\left(  R\right)
=-\sum_{j=2}^{N}\lambda_{j}\left(  R\right)  $, where the eigenvalues are
ordered as $\lambda_{N}\leq$ $\lambda_{N-1}\leq\cdots\leq\lambda_{1}$.

\subsubsection{Special cases}

We illustrate that more than one critical vector obeys $\lambda_{\max}\left(
R\right)  =1$. The particular example of the complete graph is discussed in
Section \ref{sec_critical_vector_KN}.

\textbf{1.} The homogeneous threshold $\tau_{\text{hom;}c}$ is found when
$\tau_{i}=\tau$, in which case $\lambda_{\max}\left(  R\right)  =1$ reduces to
$\frac{1}{\tau_{\text{hom;}c}}=\lambda_{\max}\left(  A\right)  $, a basic
result in \cite{PVM_ToN_VirusSpread}.

\textbf{2.} When $\frac{\delta_{i}}{\beta_{i}}=\frac{1}{\tau_{i}}=d_{i}$ for
all $1\leq i\leq N$, we observe that $\mathcal{Q}\left(  d_{i}\right)  =Q$ if
$v_{i\infty}=\varepsilon>0$, where $\varepsilon$ is arbitrary small. In that
case, the steady-state vector is $\widetilde{V}_{\infty}\rightarrow\varepsilon
u$, thus $V_{\infty}=\varepsilon\left(  \beta_{1},\beta_{2},\ldots,\beta
_{N}\right)  $ and the critical vector $\tau_{c}=\left(  \frac{1}{d_{1}}%
,\frac{1}{d_{2}},\ldots,\frac{1}{d_{N}}\right)  $. In that case, $R=$
diag$\left(  \sqrt{\frac{1}{d_{i}}}\right)  A$diag$\left(  \sqrt{\frac
{1}{d_{i}}}\right)  $ and after a similarity transform $H_{1}=$ diag$\left(
\sqrt{\frac{1}{d_{i}}}\right)  $, we obtain the stochastic matrix \cite[p.
484-486]{PVM_PerformanceAnalysisCUP}%
\[
H_{1}RH_{1}^{-1}=\Delta^{-1}A
\]
whose largest eigenvalue is, indeed, equal to one.

\subsection{Bounding $\lambda_{\max}\left(  R\right)  $}

Applying the general Rayleigh formulation for any matrix $M$,%
\[
\lambda_{\max}=\sup_{x\neq0}\frac{x^{T}Mx}{x^{T}x}%
\]
and, knowing that all components of the eigenvector belonging to the largest
eigenvalue are non-negative, we obtain%
\[
\lambda_{\max}\left(  R\right)  =\sup_{x\neq0}\frac{x^{T}\text{diag}\left(
\sqrt{\tau_{i}}\right)  A\text{diag}\left(  \sqrt{\tau_{i}}\right)  x}{x^{T}x}%
\]
Let $z=$ diag$\left(  \sqrt{\tau_{i}}\right)  x$, then%
\begin{equation}
\lambda_{\max}\left(  R\right)  =\sup_{z\neq0}\frac{z^{T}Az}{z^{T}%
\text{diag}\left(  \frac{1}{\tau_{i}}\right)  z} \label{lamda_max(R)_Rayleigh}%
\end{equation}
If $x$ is the eigenvector of $R$ belonging to the eigenvalue $\lambda_{\max
}\left(  R\right)  =1$, then (\ref{lamda_max(R)_Rayleigh}) implies that the
vector $z$ satisfies%
\[
z^{T}\text{diag}\left(  \frac{1}{\tau_{i}}\right)  z=z^{T}Az
\]
which shows that $z$ (with positive vector components) cannot be an
eigenvector of $A$, unless all $\tau_{i}=\tau$. Indeed, suppose that $z$ is an
eigenvector of $A$ belonging to $\lambda\left(  A\right)  $, then
$z^{T}Az=\lambda\left(  A\right)  z^{T}z$, which can only be equal to $z^{T}%
$diag$\left(  \frac{1}{\tau_{i}}\right)  z$ if all $\tau_{i}=\tau$ and
$\lambda\left(  A\right)  =\lambda_{\max}\left(  A\right)  =\frac{1}{\tau}$;
thus, only in the homogeneous case. In the sequel, we deduce several bounds
from (\ref{lamda_max(R)_Rayleigh}).

First, we rewrite (\ref{lamda_max(R)_Rayleigh}) as%
\begin{align*}
\lambda_{\max}\left(  R\right)   &  =\sup_{z\neq0}\frac{z^{T}Az}{z^{T}z}%
\frac{z^{T}z}{z^{T}\text{diag}\left(  \frac{1}{\tau_{i}}\right)  z}\\
&  \geq\sup_{z\neq0}\frac{z^{T}Az}{z^{T}z}\sup_{z\neq0}\frac{z^{T}z}%
{z^{T}\text{diag}\left(  \frac{1}{\tau_{i}}\right)  z}\\
&  =\lambda_{\max}\left(  A\right)  \min_{1\leq j\leq N}\tau_{i}%
\end{align*}
Thus,%
\begin{equation}
\lambda_{\max}\left(  A\right)  \min_{1\leq j\leq N}\tau_{i}\leq\lambda_{\max
}\left(  R\right)  \leq\lambda_{\max}\left(  A\right)  \max_{1\leq j\leq
N}\tau_{i} \label{lower_upper_bound_lamda_max_R}%
\end{equation}
where the upper bound follows similarly from $\sup_{z\neq0}\frac{z^{T}%
Az}{z^{T}\text{diag}\left(  \frac{1}{\tau_{i}}\right)  z}\leq\frac{\max
_{z\neq0}z^{T}Az}{\min_{z\neq0}z^{T}\text{diag}\left(  \frac{1}{\tau_{i}%
}\right)  z}$. At the critical threshold where $\lambda_{\max}\left(
R\right)  =1$, the bounds reduce, with $\tau_{\min}=\min_{1\leq j\leq N}%
\tau_{i}$ and $\tau_{\max}=\max_{1\leq j\leq N}\tau_{i}$, to the inequality
for the minimum and maximum component of the critical $\tau$-vector,%
\[
\tau_{\min;c}\leq\frac{1}{\lambda_{\max}\left(  A\right)  }\leq\tau_{\max;c}%
\]
Hence, there is always at least one $\tau$-component below and one $\tau
$-component above the critical threshold of the homogeneous case
$\tau_{\text{hom;c}}=$ $\frac{1}{\lambda_{\max}\left(  A\right)  }$.

Next, a common lower bound (see e.g.
\cite{PVM_LAA_lowerbound_eig_symm_matrix,PVM_LAA_lowerbound_eig_symm_matrix_StephenWalker,PVM_graphspectra}%
) is obtained by letting $z=u$, the all-one vector, in
(\ref{lamda_max(R)_Rayleigh}). Equality in (\ref{lamda_max(R)_Rayleigh}) is
only achieved when $z$ is the eigenvector such that, in all other cases,%
\begin{equation}
\lambda_{\max}\left(  R\right)  \geq\frac{u^{T}Au}{u^{T}\text{diag}\left(
\frac{1}{\tau_{i}}\right)  u}=\frac{2L}{\sum_{j=1}^{N}\frac{1}{\tau_{j}}}
\label{lambda_max_R_ineq_pure_inverse_tau}%
\end{equation}
For all regular graphs\footnote{In a regular graph \cite{PVM_graphspectra},
each node has the same degree $d_{i}=d$.}, the bound
(\ref{lambda_max_R_ineq_pure_inverse_tau}) is very sharp, because $u$ is the
largest eigenvector of $A$ belonging to $\lambda_{\max}\left(  A\right)  =d$.
However, all eigenvectors of diag$\left(  \frac{1}{\tau_{i}}\right)  $ are the
basic vectors $e_{j}$ with all components equal to zero, except for the $j$-th
one that is equal to one. Written in terms of the average degree $E\left[
D\right]  =\frac{2L}{N}$ and the harmonic mean $E\left[  \tau^{-1}\right]
=\frac{1}{N}\sum_{j=1}^{N}\frac{1}{\tau_{j}}$ yields%
\[
\lambda_{\max}\left(  R\right)  \geq\frac{E\left[  D\right]  }{E\left[
\tau^{-1}\right]  }%
\]
such that at the critical threshold, where $\lambda_{\max}\left(  R\right)
=1$, there holds that $E\left[  \tau_{c}^{-1}\right]  \geq E\left[  D\right]
$. Unfortunately, the harmonic, geometric and arithmetic mean
inequality\footnote{For real positive numbers $a_{1},a_{2,},\ldots,a_{n}$, the
harmonic, geometric and arithmetic mean inequality is
\begin{equation}
\frac{n}{\sum_{j=1}^{n}\frac{1}{a_{j}}}\leq\sqrt[n]{%
{\displaystyle\prod\limits_{j=1}^{n}}
a_{j}}\leq\frac{1}{n}\sum_{j=1}^{n}a_{j}
\label{harmonic_geometric_arithmetic_mean_inequality}%
\end{equation}
}, that leads to $\frac{1}{E\left[  \tau^{-1}\right]  }=N\left(  \sum
_{j=1}^{N}\frac{1}{\tau_{j}}\right)  ^{-1}\leq\frac{1}{N}\sum_{j=1}^{n}%
\tau_{j}=E\left[  \tau\right]  $, prevents us to clearly upper bound the
average zero infection $\tau$-region, $\left[  0,E\left[  \tau_{c}\right]
\right]  $. Approximative, by assuming $\frac{1}{E\left[  \tau^{-1}\right]
}\approx E\left[  \tau\right]  $, the average zero infection $\tau$-region is
upper bounded by the mean degree $E\left[  D\right]  $. Notice that, in the
homogeneous case ($\tau_{j}=\tau$), the approximation is exact, leading to
$\tau_{\text{hom;c}}\leq\frac{1}{E\left[  D\right]  }$.

There are several other interesting choices. A first alternative choice is
$z=D$, where $D=\left(  d_{1},d_{2},\ldots,d_{N}\right)  $ is the degree
vector. The Rayleigh expression (\ref{lamda_max(R)_Rayleigh}) becomes%
\[
\lambda_{\max}\left(  R\right)  \geq\frac{D^{T}AD}{D^{T}\text{diag}\left(
\frac{1}{\tau_{i}}\right)  D}=\frac{\sum_{k=1}^{N}\sum_{j=1}^{N}d_{k}%
a_{kj}d_{j}}{\sum_{j=1}^{N}\frac{d_{j}^{2}}{\tau_{j}}}%
\]
With $d_{j}=\sum_{l=1}^{N}a_{jl}$, and using symmetry, $a_{ij}=a_{ji}$,%
\begin{align*}
D^{T}AD  &  =\sum_{k=1}^{N}\sum_{j=1}^{N}\sum_{l=1}^{N}\sum_{q=1}^{N}%
a_{jl}a_{kq}a_{kj}\\
&  =\sum_{j=1}^{N}\sum_{l=1}^{N}\sum_{q=1}^{N}a_{jl}\sum_{k=1}^{N}a_{qk}%
a_{kj}\\
&  =\sum_{j=1}^{N}\sum_{l=1}^{N}\sum_{q=1}^{N}a_{jl}\left(  A^{2}\right)
_{qj}=\sum_{l=1}^{N}\sum_{q=1}^{N}\left(  A^{3}\right)  _{lq}=N_{3}%
\end{align*}
where $N_{3}$ equals the total number of walks of length $3$ in the graph.
Thus, at the critical threshold where $\lambda_{\max}\left(  R\right)  =1$,
\begin{equation}
\sum_{j=1}^{N}\frac{d_{j}^{2}}{\tau_{j}}\geq N_{3}
\label{lambda_max_R_ineq_degree_kwadraat_inverse_tau}%
\end{equation}
Invoking the Cauchy-Schwarz inequality (see e.g. \cite[p. 90]%
{PVM_PerformanceAnalysisCUP}), we further obtain%
\begin{align*}
\sum_{j=1}^{N}\frac{1}{\tau_{j}^{2}}  &  \geq\frac{N_{3}^{2}}{\sum_{j=1}%
^{N}d_{j}^{4}}\\
\sum_{j=1}^{N}\frac{1}{\tau_{j}}  &  \geq\frac{N_{3}^{2}}{\sum_{j=1}^{N}%
\frac{d_{j}^{4}}{\tau_{j}}}%
\end{align*}
A second alternative choice is to choose the components of the vector $z$
equal to a row vector of $A$, i.e. $z_{j}=a_{qj}$, such that%
\[
\lambda_{\max}\left(  R\right)  \geq\frac{\sum_{k=1}^{N}\sum_{j=1}^{N}%
a_{qk}a_{kj}a_{qj}}{\sum_{j=1}^{N}\frac{a_{qj}^{2}}{\tau_{j}}}%
\]
Since%
\[
\sum_{k=1}^{N}\sum_{j=1}^{N}a_{qk}a_{kj}a_{qj}=\sum_{j=1}^{N}\left(
A^{2}\right)  _{qj}a_{qj}=\left(  A^{3}\right)  _{qq}%
\]
and $\sum_{j=1}^{N}\frac{a_{qj}^{2}}{\tau_{j}}=\sum_{j=1}^{N}\frac{a_{qj}%
}{\tau_{j}}$, we obtain at the critical threshold where $\lambda_{\max}\left(
R\right)  =1$,%
\[
\sum_{j=1}^{N}\frac{a_{qj}}{\tau_{j}}\geq\left(  A^{3}\right)  _{qq}%
\]
Summing over all $q$ leads to%
\begin{equation}
\sum_{j=1}^{N}\frac{d_{j}}{\tau_{j}}\geq\text{ trace}\left(  A^{3}\right)
\label{lambda_max_R_ineq_degree_inverse_tau}%
\end{equation}

\subsection{Computation of $\lambda_{\max}\left(  R\right)  $ in $K_{N}$}

\label{sec_critical_vector_KN}The adjacency matrix of the complete graph
$K_{N}$ is $A_{K_{N}}=J-I$, where $J=u.u^{T}$ is the all-one matrix. Then, the
$R$ matrix defined in (\ref{def_R}), is%
\begin{align*}
R_{K_{N}}  &  =\text{diag}\left(  \sqrt{\tau_{i}}\right)  \left(  J-I\right)
\text{diag}\left(  \sqrt{\tau_{i}}\right) \\
&  =\text{diag}\left(  \sqrt{\tau_{i}}\right)  u.u^{T}\text{diag}\left(
\sqrt{\tau_{i}}\right)  -\text{diag}\left(  \tau_{i}\right) \\
&  =\left(  u^{T}\text{diag}\left(  \sqrt{\tau_{i}}\right)  \right)
^{T}.u^{T}\text{diag}\left(  \sqrt{\tau_{i}}\right)  -\text{diag}\left(
\tau_{i}\right) \\
&  =\sqrt{\tau}.\sqrt{\tau}^{T}-\text{diag}\left(  \tau_{i}\right)
\end{align*}
where the square root vector of $\tau$ is $\sqrt{\tau}=\left(  \sqrt{\tau_{1}%
},\sqrt{\tau_{2}},\ldots,\sqrt{\tau_{N}}\right)  $. The eigenvalues are
determined by the zeros of the characteristic polynomial $p_{N}\left(
\lambda\right)  =\det\left(  R-\lambda I\right)  $,%
\begin{align*}
p_{N}\left(  \lambda\right)   &  =\det\left(  \sqrt{\tau}.\sqrt{\tau}%
^{T}-\text{diag}\left(  \tau_{i}+\lambda\right)  \right) \\
&  =\det\left(  -\text{diag}\left(  \tau_{i}+\lambda\right)  \right) \\
&  \hspace{0.5cm}\times\det\left(  I-\text{diag}\left(  \frac{1}{\tau
_{i}+\lambda}\right)  \sqrt{\tau}.\sqrt{\tau}^{T}\right)
\end{align*}
After using the one-rank update formula (see e.g. \cite{Meyer_matrix}),
$\det\left(  I+cd^{T}\right)  =1+d^{T}c$, we obtain%
\begin{align*}
p_{N}\left(  \lambda\right)   &  =(-1)^{N}\left(  1-\sqrt{\tau}^{T}%
.\text{diag}\left(  \frac{1}{\tau_{i}+\lambda}\right)  \sqrt{\tau}\right)
{\displaystyle\prod\limits_{i=1}^{N}}
\left(  \tau_{i}+\lambda\right) \\
&  =(-1)^{N}\left(  1-\sum_{j=1}^{N}\frac{\tau_{j}}{\tau_{j}+\lambda}\right)
{\displaystyle\prod\limits_{i=1}^{N}}
\left(  \tau_{i}+\lambda\right)
\end{align*}
Let us order the non-negative vector components of $\tau$ as $0\leq
\tau_{\left(  N\right)  }\leq\tau_{\left(  N-1\right)  }\leq\cdots\leq
\tau_{\left(  1\right)  }$. The rational function $r\left(  \lambda\right)
=1-\sum_{j=1}^{N}\frac{\tau_{j}}{\tau_{j}+\lambda}$ has simple poles at
$\lambda=-\tau_{j}$ and is increasing between two consecutive poles. Moreover,
$\lim_{\lambda\rightarrow\pm\infty}r\left(  \lambda\right)  =1$. This implies
that $r\left(  \lambda\right)  $ has simple zeros between each pair $\left(
-\tau_{\left(  j-1\right)  },-\tau_{\left(  j\right)  }\right)  $ and those
zeros are the zeros of the characteristic polynomial $p_{N}\left(
\lambda\right)  =r\left(  \lambda\right)
{\displaystyle\prod\limits_{i=1}^{N}}
\left(  \tau_{i}+\lambda\right)  $ provided $p_{N}\left(  -\tau_{j}\right)
=-\tau_{j}%
{\displaystyle\prod\limits_{i=1;i\neq j}^{N}}
\left(  \tau_{i}-\tau_{j}\right)  \neq0$, i.e. provided all $\tau_{i}$ are
different. The largest zero of $p_{N}\left(  \lambda\right)  $ exceeds
$\lambda=-\tau_{\left(  N\right)  }\leq0$. Even much sharper, since
trace$\left(  A\right)  =\sum_{i=1}^{N}\lambda_{i}=0$, we know that%
\[
\lambda_{\max}=\lambda_{1}=-\sum_{i=1}^{N-1}\lambda_{i}\leq\sum_{i=1}^{N}%
\tau_{i}-\tau_{\min}%
\]

We rewrite $r\left(  \lambda\right)  $ as
\[
r\left(  \lambda\right)  =\lambda\sum_{j=1}^{N}\frac{1}{\tau_{j}+\lambda
}-\left(  N-1\right)
\]
from which the largest zero of $(-1)^{N}p_{N}\left(  \lambda\right)  $ is the
only positive solution in $\lambda$ of%
\begin{equation}
\sum_{j=1}^{N}\frac{1}{\tau_{j}+\lambda}=\frac{N-1}{\lambda}
\label{equation_largest_eig_R_KN}%
\end{equation}
By iteration of the rewritten equation as $\lambda=\frac{1}{\frac{1}{N-1}%
\sum_{j=1}^{N}\frac{1}{\tau_{j}+\lambda}}$, we obtain the continued fraction%
\[
\lambda_{\max}=\frac{1}{\frac{1}{N-1}\sum_{j=1}^{N}\frac{1}{\tau_{j}+\frac
{1}{\frac{1}{N-1}\sum_{k=1}^{N}\frac{1}{\tau_{k}+\ddots\frac{\ddots}{\tau
_{q}+\frac{1}{\frac{1}{N-1}\sum_{l=1}^{N}\frac{1}{\tau_{l}+\ddots}}}}}}}%
\]
from which the following convergents are deduced,%
\[
\frac{N-1}{\sum_{j=1}^{N}\frac{1}{\tau_{j}}}<\frac{N-1}{\sum_{j=1}^{N}\frac
{1}{\tau_{j}+\frac{N-1}{\sum_{k=1}^{N}\frac{1}{\tau_{k}}}}}<\cdots\leq
\lambda_{\max}%
\]
Notice that these convergents for $K_{N}$ show that, indeed,
(\ref{lambda_max_R_ineq_pure_inverse_tau}) is a sharp bound for regular
graphs. Lagrange expansion of (\ref{equation_largest_eig_R_KN}) is also
possible, but we omit this analysis.

The critical vector components thus satisfy, with $\lambda_{\max}\left(
R\right)  =1$, the equation%
\begin{equation}
\sum_{j=1}^{N}\frac{1}{\tau_{j}+1}=N-1
\label{critical_equation_complete_graph}%
\end{equation}
A critical $\tau$-vector must have bounded components. For, if $\tau
_{k}\rightarrow\infty$, then (\ref{critical_equation_complete_graph}) implies
that all other $\tau_{j}=0$, which leads to a physically uninteresting
situation. Let $\tau_{j}=$ $\tau_{\hom;c}+h_{j}$, where $\tau_{\hom;c}%
=\frac{1}{N-1}$ as shown below, then (\ref{critical_equation_complete_graph})
can be rewritten as%
\[
\sum_{j=1}^{N}\frac{1}{1+\frac{N-1}{N}h_{j}}=N
\]
For small $h_{j}$ where $\left(  1+\frac{N-1}{N}h_{j}\right)  ^{-1}%
=1-\frac{N-1}{N}h_{j}+O\left(  h_{j}^{2}\right)  $, we have that $\sum
_{j=1}^{N}h_{j}\approx0$. Hence, the small deviations $h_{j}$ from the
homogeneous case are balanced, in the sense that the net or average deviation
is about zero. Suppose that all $h_{j}=0$ for $3\leq j\leq N$, then $h_{1}$
and $h_{2}$ obey a hyperbolic relation%
\[
h_{1}=\frac{-h_{2}}{1+2\frac{N-1}{N}h_{2}}%
\]
Small negative values for $h_{2}$ correspond, on the critical threshold, to
large positive values for $h_{1}$ (and vice versa).

Finally, the homogeneous case, where $\tau_{j}=\tau_{\hom}$, considerably
simplifies to the characteristic polynomial%
\[
p_{N}\left(  \lambda\right)  =(-1)^{N}\left(  \lambda-\tau_{\hom}\left(
N-1\right)  \right)  \left(  \tau_{\hom}+\lambda\right)  ^{N-1}%
\]
whose zeros are $\lambda=\tau_{\hom}\left(  N-1\right)  $ and $\lambda
=-\tau_{\hom}$ with multiplicity $N-1$. This example illustrates that,
although heterogeneity is much more natural, it complicates analysis seriously.

\subsection{Additional properties}

We list here additional properties that have been proved in
\cite{PVM_ToN_VirusSpread}, and whose extension to the in-homogenous setting
is rather straightforward.

\begin{lemma}
\label{lemma_all_vi_zero_or_none_are_zero}In a connected graph, either
$v_{i\infty}=0$ for all $i$ nodes, or none of the components $v_{i\infty}$ is zero.
\end{lemma}

Lemma \ref{lemma_all_vi_zero_or_none_are_zero} also follows from the
Perron-Frobenius theorem as shown in the proof of Theorem
\ref{theorem_all_eigenvalues_generalized_Q_are_positive_except_1}.

\begin{theorem}
\label{Theorem_continued_fraction}The non-zero steady-state infection
probability of any node $i$ in the $N$-intertwined model can be expressed as a
continued fraction{\footnotesize
\begin{equation}
v_{i\infty}=1-\frac{1}{1+\frac{\gamma_{i}}{\delta_{i}}-\delta_{i}^{-1}%
\sum_{j=1}^{N}\frac{\beta_{j}a_{ij}}{1+\frac{\gamma_{j}}{\delta_{j}}%
-\delta_{j}^{-1}\sum_{k=1}^{N}\frac{\beta_{k}a_{jk}}{1+\frac{\gamma_{k}%
}{\delta_{k}}-\delta_{k}^{-1}\sum_{q=1}^{N}\frac{a_{qk}\beta_{q}}{\ddots}}}}
\label{v_infinite_fraction}%
\end{equation}
} where the total infection rate of node $i$, incurred by all neighbors
towards node $i$, is%
\begin{equation}
\gamma_{i}=\sum_{j=1}^{N}a_{ij}\beta_{j}=\sum_{j\in\text{ neighbor}\left(
i\right)  }\beta_{j} \label{def_gamma_i}%
\end{equation}
Consequently, the exact steady-state infection probability of any node $i$ is
bounded by%
\begin{equation}
0\leq v_{i\infty}\leq1-\frac{1}{1+\frac{\gamma_{i}}{\delta_{i}}}
\label{bounds_on_vi_inf}%
\end{equation}

\end{theorem}

As explained in \cite{PVM_ToN_VirusSpread}, the continued fraction stopped at
iteration $k$ includes the effect of virus spread up to the ($k-1)$-hop
neighbors of node $i$. In the homogeneous case where $\beta_{j}=\beta$ for all
$1\leq j\leq N$, we have that $\gamma_{i}=\beta d_{i}$ is proportional to the
degree of node \thinspace$i$. The ratio $\widetilde{\tau}_{i}=\frac{\gamma
_{i}}{\delta_{i}}$ is the total effective infection rate of node $i$.

\begin{lemma}
\label{lem_lower_bound} In a connected graph $G$ above the critical threshold,
a lower bound of $v_{i\infty}$ for any node $i$ equals%
\begin{equation}
v_{i\infty}\geq1-\frac{1}{\min_{1\leq k\leq N}\frac{\gamma_{k}}{\delta_{k}}}
\label{lower_bound_any_vi_infty}%
\end{equation}

\end{lemma}

\textbf{Proof:} Lemma \ref{lemma_all_vi_zero_or_none_are_zero} and Theorem
\ref{theorem_critical_threshold} show that, for vectors $\tau$ above the
critical threshold vector $\tau_{c}$, there exists a non-zero minimum
$v_{\min}=\min_{1\leq i\leq N}v_{i\infty}>0$ of the steady-state infection
probabilities, which obeys (\ref{steady_state_node_i_equation}). Assuming that
this minimum $v_{\min}$ occurs at node $i$,%
\[
v_{\min}=1-\frac{1}{1+\delta_{i}^{-1}\sum_{j=1}^{N}a_{ij}\beta_{j}v_{j\infty}%
}\geq1-\frac{1}{1+\frac{\gamma_{i}}{\delta_{i}}v_{\min}}%
\]
where we have used the definition (\ref{def_gamma_i}). From the last
inequality, it follows that%
\begin{equation}
v_{\min}\geq1-\frac{\delta_{i}}{\gamma_{i}} \label{bound_for_v_min}%
\end{equation}
such that (\ref{lower_bound_any_vi_infty}) is proved.\hspace{1cm}$\square$

By combining (\ref{bounds_on_vi_inf}) and (\ref{lower_bound_any_vi_infty}),
the total fraction of infected nodes $y_{\infty}=\frac{1}{N}\sum_{k=1}%
^{N}v_{k\infty}$ in steady-state is bounded by%
\[
1-\frac{1}{\min_{1\leq k\leq N}\frac{\gamma_{k}}{\delta_{k}}}\leq y_{\infty
}\leq1-\frac{1}{N}\sum_{i=1}^{N}\frac{1}{1+\frac{\gamma_{i}}{\delta_{i}}}%
\]

\section{The convexity of $v_{i\infty}$ as a function of $\delta_{i}$}

\label{sec_convexity}It is of interest (e.g. in game theory
\cite{PVM_Jasmina_Game_protection_INFOCOM2009}) to know whether the
steady-state infection probability $v_{i\infty}$ is convex in the own curing
rate $\delta_{i}$, given that all other curing rates $\delta_{j}$ for $1\leq
j\neq i\leq N$ are constant. In many infection situations, the node $i$ cannot
control the spreading process, but it can protect itself better by increasing
its own curing rate $\delta_{i}$, for example, by installing more effective
antivirus software in computer networks, or by vaccinating people against some diseases.

\begin{theorem}
\label{theorem_convexity}If all curing rates are the same, i.e. $\delta
_{k}=\delta_{i}$ for $1\leq k\leq N$, then $v_{k\infty}$ is convex in
$\delta_{i}$.
\end{theorem}

However, if all curing rates $\delta_{j}$ for $1\leq j\neq i\leq N$ are
constant and independent from each other and from the infection rates
$\beta_{j}$, the non-zero steady-state infection probability $v_{k\infty
}\left(  \delta_{1},\ldots,\delta_{i},\ldots,\delta_{N}\right)  >0$ \emph{can
be} concave in $\delta_{i}$.

\textbf{Proof:} We operate above the critical threshold specified by
$\lambda_{\max}\left(  R\right)  =1$, where the vector $V_{\infty}>0$ and
start from the steady-state equation (\ref{steady_state_node_i_equation}) for
node $i$. Differentiation with respect to $\delta_{i}$ results in%
\begin{equation}
\sum_{k=1}^{N}a_{ik}\beta_{k}\frac{\partial v_{k\infty}}{\partial\delta_{i}%
}=\frac{v_{i\infty}}{1-v_{i\infty}}+\frac{\delta_{i}}{\left(  1-v_{i\infty
}\right)  ^{2}}\frac{\partial v_{i\infty}}{\partial\delta_{i}}
\label{eerste_afg_naar_delta_i}%
\end{equation}
and%
\begin{align*}
\sum_{k=1}^{N}a_{ik}\beta_{k}\frac{\partial^{2}v_{k\infty}}{\partial\delta
_{i}^{2}}  &  =\frac{2}{\left(  1-v_{i\infty}\right)  ^{2}}\frac{\partial
v_{i\infty}}{\partial\delta_{i}}+\frac{2\delta_{i}}{\left(  1-v_{i\infty
}\right)  ^{3}}\left(  \frac{\partial v_{i\infty}}{\partial\delta_{i}}\right)
^{2}\\
&  \hspace{0.5cm}+\frac{\delta_{i}}{\left(  1-v_{i\infty}\right)  ^{2}}%
\frac{\partial^{2}v_{i\infty}}{\partial\delta_{i}^{2}}%
\end{align*}
Differentiating any other row $j\neq i$ in (\ref{steady_state_equation})
\[
\sum_{k=1}^{N}a_{jk}\beta_{k}v_{k\infty}=\frac{v_{j\infty}}{1-v_{j\infty}%
}\delta_{j}%
\]
with respect to $\delta_{i}$ results in%
\[
\sum_{k=1}^{N}a_{jk}\beta_{k}\frac{\partial v_{k\infty}}{\partial\delta_{i}%
}=\frac{\delta_{j}}{\left(  1-v_{j\infty}\right)  ^{2}}\frac{\partial
v_{j\infty}}{\partial\delta_{i}}%
\]
and%
\[
\sum_{k=1}^{N}a_{jk}\beta_{k}\frac{\partial^{2}v_{k\infty}}{\partial\delta
_{i}^{2}}=\frac{2\delta_{j}\left(  \frac{\partial v_{j\infty}}{\partial
\delta_{i}}\right)  ^{2}}{\left(  1-v_{j\infty}\right)  ^{3}}+\frac{\delta
_{j}}{\left(  1-v_{j\infty}\right)  ^{2}}\frac{\partial^{2}v_{j\infty}%
}{\partial\delta_{i}^{2}}%
\]

Written in matrix form, we have%
\begin{equation}
A\text{diag}\left(  \beta_{k}\right)  \frac{\partial V_{\infty}}%
{\partial\delta_{i}}=\text{diag}\left(  \frac{\delta_{k}}{\left(
1-v_{k\infty}\right)  ^{2}}\right)  \frac{\partial V_{\infty}}{\partial
\delta_{i}}+\frac{v_{i\infty}}{1-v_{i\infty}}e_{i}
\label{matrix_eerste_afg_naar_delta_i}%
\end{equation}
where the basisvector $e_{i}$ has all zero components, except for the
component $i$ that equals 1. When curing rate $\delta_{j}$ is a function of
$\delta_{k}$, the equations change. In particular, if $\delta_{k}=\delta_{i}$
for all $1\leq k\leq N$, then the vector $e_{i}$ must be replaced by the
all-one vector $u$.

The second order derivatives are, in matrix form,%
\begin{align*}
A\text{diag}\left(  \beta_{k}\right)  \frac{\partial^{2}V_{\infty}}%
{\partial\delta_{i}^{2}}  &  =W_{i\infty}+\text{diag}\left(  \frac{\delta_{k}%
}{\left(  1-v_{k\infty}\right)  ^{2}}\right)  \frac{\partial^{2}V_{\infty}%
}{\partial\delta_{i}^{2}}\\
&  \hspace{0.5cm}+\frac{2}{\left(  1-v_{i\infty}\right)  ^{2}}\frac{\partial
v_{i\infty}}{\partial\delta_{i}}e_{i}%
\end{align*}
where $W_{i\infty}=\left[
\begin{array}
[c]{ccc}%
\frac{2\delta_{1}\left(  \frac{\partial v_{1\infty}}{\partial\delta_{i}%
}\right)  ^{2}}{\left(  1-v_{1\infty}\right)  ^{3}} & \cdots & \frac
{2\delta_{N}\left(  \frac{\partial v_{N\infty}}{\partial\delta_{i}}\right)
^{2}}{\left(  1-v_{N\infty}\right)  ^{3}}%
\end{array}
\right]  ^{T}$.

We rewrite the matrix equations as%
\begin{equation}
S\frac{\partial V_{\infty}}{\partial\delta_{i}}=-\frac{v_{i\infty}%
}{1-v_{i\infty}}e_{i} \label{matrix_eq_derivatives_v_naar_delta}%
\end{equation}
where the matrix%
\begin{equation}
S=\text{diag}\left(  \frac{\delta_{j}}{\left(  1-v_{j\infty}\right)  ^{2}%
}\right)  -A\text{diag}\left(  \beta_{k}\right)  \label{def_S}%
\end{equation}
is written in terms of the generalized Laplacian $\mathcal{Q}\left(
q_{i}\right)  $, defined in (\ref{def_generalized_Laplacian}), as\footnote{We
remark that, with $B=$diag$\left(  \sqrt{\beta_{k}}\right)  $, the matrix
\[
BSB^{-1}=\text{diag}\left(  \sqrt{\beta_{k}}\right)  \mathcal{Q}\left(
\frac{1}{\tau_{j}\left(  1-v_{j\infty}\right)  ^{2}}\right)  \text{diag}%
\left(  \sqrt{\beta_{k}}\right)
\]
is symmetric.}%
\begin{equation}
S=\mathcal{Q}\left(  \frac{1}{\tau_{j}\left(  1-v_{j\infty}\right)  ^{2}%
}\right)  \text{diag}\left(  \beta_{j}\right)  \label{T_in_terms_of_Q}%
\end{equation}
Lemma \ref{lemma_generalized_laplacian_q*>q_positive_def} shows that $S$ is
positive definite, which implies that also $S^{-1}$ is positive definite
because $S=U$diag$\left(  \lambda_{j}\right)  U^{T}$ shows that $S^{-1}%
=U$diag$\left(  \lambda_{j}^{-1}\right)  U^{T}$ and, thus, that the inverse
$S^{-1}$ exists. The vector $\frac{\partial V_{\infty}}{\partial\delta_{i}}$
is solved from (\ref{matrix_eq_derivatives_v_naar_delta}) explicitly as
\begin{equation}
\frac{\partial V_{\infty}}{\partial\delta_{i}}=-\frac{v_{i\infty}%
}{1-v_{i\infty}}S^{-1}e_{i}=-\frac{v_{i\infty}}{1-v_{i\infty}}\left(
S^{-1}\right)  _{\text{column }i} \label{solution_DVectorV_delta_i}%
\end{equation}
from which%
\begin{equation}
\left(  S^{-1}\right)  _{ki}=\left(  1-\frac{1}{v_{i\infty}}\right)
\frac{\partial v_{k\infty}}{\partial\delta_{i}} \label{element_(inverseS)}%
\end{equation}
Increasing the virus curing rate cannot increase the virus infection
probability, such that $\frac{\partial v_{k\infty}}{\partial\delta_{i}}\leq0$
for all $1\leq k\leq N$. This implies that all elements of $S^{-1}$ are
non-negative. Moreover, since $\frac{\partial v_{k\infty}}{\partial\delta_{i}%
}\leq0$, the left-hand side in (\ref{eerste_afg_naar_delta_i}) is always
negative, which leads, in a different way, to the inequality
(\ref{inequality_inverse}).

Only at the critical threshold, the derivatives $\frac{\partial V_{\infty}%
}{\partial\delta_{i}}$ do not exist because the left- and right derivative at
that point are not equal. Below the critical threshold, where $V_{\infty}=0$,
(\ref{solution_DVectorV_delta_i}) does not yield information about the
existence of $S^{-1}$. However, the definition (\ref{def_S}) shows that $S=$
diag$\left(  \delta_{j}\right)  -A$diag$\left(  \beta_{j}\right)  $. Hence, if
$\frac{\delta_{j}}{\beta_{j}}=d_{j}$ for each node $j$, then diag$\left(
\beta_{j}^{-1}\right)  S$ equals the Laplacian $Q$ and $S^{-1}$ does not
exist. In general, it is difficult to conclude for which vector $C=\left(
\delta_{1},\delta_{2},\ldots,\delta_{N}\right)  $ that $S^{-1}$ exists below
the critical threshold. But, below the critical threshold, $V_{\infty}=0$ such
that both convexity and concavity hold. In the sequel, we ignore further
considerations about this sub-threshold regime.

We recast the second order derivatives\footnote{In fact, we can show, for any
integer $m>0$, that%
\[
S\frac{\partial^{m}V_{\infty}}{\partial\delta_{i}^{m}}=R_{m}%
\]
so that any higher order derivative vector equals%
\[
\frac{\partial^{m}V_{\infty}}{\partial\delta_{i}^{m}}=S^{-1}R_{m}%
\]
which illustrates the importance of the positive definite matrix $S$ and its
non-negative inverse $S^{-1}$.} in terms of the matrix $S$,%
\[
S\frac{\partial^{2}V_{\infty}}{\partial\delta_{i}^{2}}=-\widetilde{W}%
\]
where%
\[
\widetilde{W}=W_{i\infty}+\frac{2}{\left(  1-v_{i\infty}\right)  ^{2}}%
\frac{\partial v_{i\infty}}{\partial\delta_{i}}e_{i}%
\]
Above the critical threshold, $S^{-1}$ exists such that%
\begin{equation}
\frac{\partial^{2}V_{\infty}}{\partial\delta_{i}^{2}}=-S^{-1}\widetilde{W}
\label{solution_D^2_VectorV_delta_i}%
\end{equation}

Introducing (\ref{solution_DVectorV_delta_i}) in $W_{\infty}$ yields%
\[
W_{i\infty}=\frac{2v_{i\infty}^{2}}{\left(  1-v_{i\infty}\right)  ^{2}}\left[
\begin{array}
[c]{ccc}%
\frac{\delta_{1}\left(  \left(  S^{-1}\right)  _{1i}\right)  ^{2}}{\left(
1-v_{1\infty}\right)  ^{3}} & \cdots & \frac{\delta_{N}\left(  \left(
S^{-1}\right)  _{Ni}\right)  ^{2}}{\left(  1-v_{N\infty}\right)  ^{3}}%
\end{array}
\right]  ^{T}%
\]
and%
\[
\frac{2}{\left(  1-v_{i\infty}\right)  ^{2}}\frac{\partial v_{i\infty}%
}{\partial\delta_{i}}e_{i}=-\frac{2v_{i\infty}}{\left(  1-v_{i\infty}\right)
^{3}}\left(  S^{-1}\right)  _{ii}e_{i}%
\]
which shows that the right-hand side vector $\widetilde{W}$ has all positive
elements, except for the $i$-th component which is%
\begin{align*}
\widetilde{W}_{i}  &  =\frac{2\delta_{i}\left(  v_{i\infty}\left(
S^{-1}\right)  _{ii}\right)  ^{2}}{\left(  1-v_{i\infty}\right)  ^{5}}%
-\frac{2v_{i\infty}\left(  S^{-1}\right)  _{ii}}{\left(  1-v_{i\infty}\right)
^{3}}\\
&  =2\frac{v_{i\infty}\left(  S^{-1}\right)  _{ii}}{\left(  1-v_{i\infty
}\right)  ^{3}}\left\{  \delta_{i}\frac{v_{i\infty}\left(  S^{-1}\right)
_{ii}}{\left(  1-v_{i\infty}\right)  ^{2}}-1\right\}
\end{align*}
Since all elements of $S^{-1}$ are positive and Lemma
\ref{lemma_all_vi_zero_or_none_are_zero} states that all $v_{i\infty}>0$ above
the critical threshold, we conclude from
(\ref{steady_state_equation_in_terms_of_S_row_i}), derived in Appendix
\ref{sec_S_steady_state}, that
\begin{equation}
\frac{\delta_{i}v_{i\infty}\left(  S^{-1}\right)  _{ii}}{\left(  1-v_{i\infty
}\right)  ^{2}}<1 \label{inequality_steady_state}%
\end{equation}
Hence, $\widetilde{W}_{i}<0$, but $\widetilde{W}_{k}>0$ for $k\neq i$.

When all curing rates are the same (i.e. $\delta_{k}=\delta_{i}$ for all
$1\leq k\leq N$), then, as mentioned before, $e_{i}$ needs to be replaced by
$u$, so that all components of $\widetilde{W}$ are negative. Consequently,
when all curing rates are the same and equal to $\delta_{i}$, we conclude from
(\ref{solution_D^2_VectorV_delta_i}) that the steady-state infection
probability $v_{k\infty}$ (each node $k$) is convex in $\delta_{i}$. This
proves Theorem \ref{theorem_convexity}.$\hfill\square$

When all curing rates are independent from each other, the $k$-th component in
(\ref{solution_D^2_VectorV_delta_i}) equals%
\begin{align}
\frac{\partial^{2}v_{k\infty}}{\partial\delta_{i}^{2}}  &  =-\sum_{j=1}%
^{N}\left(  S^{-1}\right)  _{kj}\widetilde{W_{j}}\nonumber\\
&  =-\frac{2v_{i\infty}^{2}}{\left(  1-v_{i\infty}\right)  ^{2}}%
\sum_{j=1;j\neq i}^{N}\left(  S^{-1}\right)  _{kj}\frac{\delta_{j}\left(
\left(  S^{-1}\right)  _{ji}\right)  ^{2}}{\left(  1-v_{j\infty}\right)  ^{3}%
}\nonumber\\
&  \hspace{0.5cm}+2\frac{v_{i\infty}\left(  S^{-1}\right)  _{ki}\left(
S^{-1}\right)  _{ii}}{\left(  1-v_{i\infty}\right)  ^{3}}\left\{  1-\delta
_{i}\frac{v_{i\infty}\left(  S^{-1}\right)  _{ii}}{\left(  1-v_{i\infty
}\right)  ^{2}}\right\}  \label{tweede_afg_vk_delta_i}%
\end{align}
Hence,%
\[
\frac{\left(  1-v_{i\infty}\right)  ^{2}}{2v_{i\infty}}\frac{\partial
^{2}v_{k\infty}}{\partial\delta_{i}^{2}}=M_{ki}%
\]
where%
\begin{equation}
M_{ki}=\frac{\left(  S^{-1}\right)  _{ki}\left(  S^{-1}\right)  _{ii}}{\left(
1-v_{i\infty}\right)  }-v_{i\infty}\sum_{j=1}^{N}\left(  S^{-1}\right)
_{kj}\frac{\delta_{j}\left(  \left(  S^{-1}\right)  _{ji}\right)  ^{2}%
}{\left(  1-v_{j\infty}\right)  ^{3}} \label{def_M_ki}%
\end{equation}
Unfortunately, it is difficult in general to determine the sign of $M_{ki}$ as
further illustrated in Appendix \ref{sec_analysis_M_ki}.

Simulations show that $v_{k\infty}\left(  \delta_{1},\ldots,\delta_{i}%
,\ldots,\delta_{N}\right)  $ (for any $k$) can be convex in $\delta_{i}$ (e.g.
in the lattice and complete graph) as well as concave (e.g. in a star). These
simulations indicate that either regime is possible, but no combination (i.e.
$v_{k\infty}\left(  \delta_{1},\ldots,\delta_{i},\ldots,\delta_{N}\right)  $
is convex in some $\delta_{i}$ region, but concave in another) was encountered.

\section{The derivatives $\frac{\partial v_{i\infty}}{\partial\delta_{i}}$}

Our starting point is the matrix equation
(\ref{matrix_eq_derivatives_v_naar_delta}), which we solve here by using
Cramer's rule,%
\[
\frac{\partial v_{i\infty}}{\partial\delta_{i}}=-\frac{v_{i\infty}%
}{1-v_{i\infty}}\frac{\det\left(  S_{G\backslash\left\{  i\right\}  }\right)
}{\det S}%
\]
where $G\backslash\left\{  i\right\}  $ denotes the graph $G$ from which the
node $i$ is removed (together with all its incident links). Using the
definition (\ref{T_in_terms_of_Q}) of $S$ shows that%
\[
\frac{\partial v_{i\infty}}{\partial\delta_{i}}=-\frac{v_{i\infty}}{\beta
_{i}\left(  1-v_{i\infty}\right)  }\frac{\det\left(  \mathcal{Q}%
_{G\backslash\left\{  i\right\}  }\left(  \frac{1}{\tau_{j}\left(
1-v_{j\infty}\right)  ^{2}}\right)  \right)  }{\det\mathcal{Q}\left(  \frac
{1}{\tau_{j}\left(  1-v_{j\infty}\right)  ^{2}}\right)  }%
\]
A determinant is unchanged by interchanging two rows and two columns. This
means that we can write the matrix%
\[
\mathcal{Q}\left(  \frac{1}{\tau_{j}\left(  1-v_{j\infty}\right)  ^{2}%
}\right)  =\left[
\begin{array}
[c]{cc}%
\mathcal{Q}_{G\backslash\left\{  i\right\}  }\left(  \frac{1}{\tau_{j}\left(
1-v_{j\infty}\right)  ^{2}}\right)  & -a_{i}\\
-a_{i}^{T} & \frac{1}{\tau_{i}\left(  1-v_{i\infty}\right)  ^{2}}%
\end{array}
\right]
\]
where the vector $a_{i}$ is the relabeled connection vector of node $i$ to all
other nodes in $G$ and $a_{i}^{T}a_{i}=d_{i}$. Invoking
\begin{equation}
\det\left[
\begin{array}
[c]{cc}%
A & B\\
C & D
\end{array}
\right]  =\det A\det\left(  D-CA^{-1}B\right)  \label{determinant_blockmatrix}%
\end{equation}
where $D-CA^{-1}B$ is called the Schur complement of $A$ (see e.g.
\cite{Meyer_matrix}), we find that%
\[
\frac{\det\mathcal{Q}}{\det\mathcal{Q}_{G\backslash\left\{  i\right\}  }%
}=\frac{1}{\tau_{i}\left(  1-v_{i\infty}\right)  ^{2}}-f
\]
where the quadratic form is%
\[
f=a_{i}^{T}\mathcal{Q}_{G\backslash\left\{  i\right\}  }^{-1}\left(  \frac
{1}{\tau_{j}\left(  1-v_{j\infty}\right)  ^{2}}\right)  a_{i}%
\]
Whence,%
\begin{equation}
\frac{\partial v_{i\infty}}{\partial\delta_{i}}=-\frac{\left(  1-v_{i\infty
}\right)  v_{i\infty}}{\delta_{i}-\beta_{i}\left(  1-v_{i\infty}\right)
^{2}f} \label{afg_vi_delta_i}%
\end{equation}
The quadratic form $f$ does not dependent on $v_{i\infty}$. Moreover, Lemma
\ref{lemma_generalized_laplacian_q*>q_positive_def} implies that
$\mathcal{Q}_{G\backslash\left\{  i\right\}  }^{-1}\left(  \frac{1}{\tau
_{j}\left(  1-v_{j\infty}\right)  ^{2}}\right)  $ is positive definite (for
$V_{\infty}>0$). Hence, $f>0$. The fact that $\frac{\partial V_{\infty}%
}{\partial\delta_{i}}\leq0$ implies $1\geq\tau_{i}\left(  1-v_{i\infty
}\right)  ^{2}f$ and because the inequality holds for all $v_{i\infty}$, we
also have that $1\geq\tau_{i}f$.

The optimization of an utility or cost function of the type, that, for
example, appears in game theory (see
\cite{PVM_Jasmina_Game_protection_INFOCOM2009}),%
\[
J_{i}=c_{i}\delta_{i}+v_{i\infty}%
\]
where $c_{i}$ is price to protect a node $i$ against the spread of infections,
requires to compute the optimum $\frac{\partial J_{i}}{\partial\delta_{i}%
}=c_{i}+\frac{\partial v_{i\infty}}{\partial\delta_{i}}=0$ for all $1\leq
i\leq N$. With (\ref{afg_vi_delta_i}), this equation is solved explicitly as%
\[
\frac{\left(  1-v_{i\infty}\right)  v_{i\infty}}{c_{i}}+\beta_{i}\left(
1-v_{i\infty}\right)  ^{2}f=\delta_{i}^{\ast}%
\]
Thus, the optimal value of \ $\delta_{i}^{\ast}>\frac{\left(  1-v_{i\infty
}\right)  v_{i\infty}}{c_{i}}$. An exact computation of $\delta_{i}^{\ast}$ is
generally complex because $f=f\left(  \tau_{1},\ldots,\tau_{i-1},\tau
_{i+1},\ldots,\tau_{N}\right)  $ is a non-linear function that couples all the
$\tau_{j}$ (and $\delta_{i}$).

\section{Summary}

The heterogeneous $N$-intertwined virus spread model has been described and
analyzed in the steady-state. Since it applies to any network and any
combination of node infections and curing vectors, $B$ and $C$, we believe
that the heterogeneous $N$-intertwined virus spread model is useful for a wide
range of practical infection scenarios in networks, from computer viruses to
epidemics in social networks and in nature. The critical threshold regime is
investigated, bounds are presented and the metastable steady-state infection
probabilities are shown to be convex in the own curing rate provided all
curing rates are the same. When the latter is not the case, the metastable
steady-state infection probabilities can be either concave or convex.

\textbf{Acknowledgement}

We are grateful to Ariel Orda for useful discussions. We thank Zhi Xu for
pointing to an error in an earlier version and to Bo Qu for extensive
simulations to find concave $v_{k\infty}\left(  \delta_{1},\ldots,\delta
_{N}\right)  $. This research performed in 2008 was supported by Next
Generation Infrastructures (Bsik) and the EU FP7 project ResumeNet (project
No. 224619).

{\footnotesize
\bibliographystyle{plain}
\bibliography{cac,math,misc,net,pvm,qth,tel}

\begin{thebibliography}{10}

\bibitem{Biggs}
N.~Biggs.
\newblock {\em Algebraic Graph Theory}.
\newblock Cambridge University Press, Cambridge, U.K., 2nd edition, 1996.

\bibitem{Cvetkovic}
D.~M. Cvetkovi{\'c}, M.~Doob, and H.~Sachs.
\newblock {\em Spectra of Graphs, Theory and Applications}.
\newblock Johann Ambrosius Barth Verlag, Heidelberg, third edition, 1995.

\bibitem{Ganesh_2005}
A.~Ganesh, L.~Massouli{\'e}, and D.~Towsley.
\newblock The effect of network topology on the spread of epidemics.
\newblock {\em IEEE INFOCOM2005}, 2005.

\bibitem{GantmacherII}
F.~R. Gantmacher.
\newblock {\em The Theory of Matrices}, volume~II.
\newblock Chelsea Publishing Company, New York, 1959.

\bibitem{Meyer_matrix}
C.~D. Meyer.
\newblock {\em Matrix Analysis and Applied Linear Algebra}.
\newblock Society for Industrial and Applied Mathematics (SIAM), Philadelphia,
  2000.

\bibitem{Jasmina_variance2008}
J.~Omic and P.~Van~Mieghem.
\newblock The metastable state of a {SIS} model.
\newblock {\em Journal of computer and system science}, submitted 2008.

\bibitem{PVM_Jasmina_Game_protection_INFOCOM2009}
J.~Omic, P.~Van~Mieghem, and A.~Orda.
\newblock Game theory and computer viruses.
\newblock {\em IEEE Infocom2009}, 2009.

\bibitem{PVM_PerformanceAnalysisCUP}
P.~Van~Mieghem.
\newblock {\em Performance Analysis of Communications Networks and Systems}.
\newblock Cambridge University Press, Cambridge, U.K., 2006.

\bibitem{PVM_LAA_lowerbound_eig_symm_matrix}
P.~Van~Mieghem.
\newblock A new type of lower bound for the largest eigenvalue of a symmetric
  matrix.
\newblock {\em Linear Algebra and its Applications}, 427(1):119--129, November
  2007.

\bibitem{PVM_graphspectra}
P.~Van~Mieghem.
\newblock {\em Graph Spectra for Complex Networks}.
\newblock Cambridge University Press, Cambridge, U.K., 2011.

\bibitem{PVM_PAComplexNetsCUP}
P.~Van~Mieghem.
\newblock {\em Performance Analysis of Complex Networks and Systems}.
\newblock Cambridge University Press, Cambridge, U.K., 2014.

\bibitem{PVM_ToN_VirusSpread}
P.~Van~Mieghem, J.~Omic, and R.~E. Kooij.
\newblock Virus spread in networks.
\newblock {\em IEEE/ACM Transactions on Networking}, 17(1):1--14, February
  2009.

\bibitem{PVM_LAA_lowerbound_eig_symm_matrix_StephenWalker}
S.~G. Walker and P.~Van~Mieghem.
\newblock On lower bounds for the largest eigenvalue of a symmetric matrix.
\newblock {\em Linear Algebra and its Applications}, 429(2-3):519--526, July
  2008.

\bibitem{YWang2003}
Y.~Wang, D.~Chakrabarti, C.~Wang, and C.~Faloutsos.
\newblock Epidemic spreading in real networks: An eigenvalue viewpoint.
\newblock {\em 22nd International Symposium on Reliable Distributed Systems
  (SRDS'03)- IEEE Computer}, pages 25--34, October 2003.

\bibitem{Wilkinson}
J.~H. Wilkinson.
\newblock {\em The Algebraic Eigenvalue Problem}.
\newblock Oxford University Press, New York, 1965.

\end{thebibliography}
}

\appendix

\section{Properties of the matrix $S$ defined in (\ref{def_S})}

Due to the fundamental role of the positive definite matrix $S$, defined in
(\ref{def_S}), and its inverse $S^{-1}$, we present more properties.

\subsection{Deductions from the inverse of a matrix}

\label{sec_inverse_S}The $i$-th row in the identity $S^{-1}S=I$
is\footnote{Since both matrix and inverse commute, we can also consider
$SS^{-1}=I$, which leads to a slightly less simple form. Indeed,%
\[
1_{\left\{  i=j\right\}  }=\sum_{k=1}^{N}S_{ik}\left(  S^{-1}\right)
_{kj}=S_{ii}\left(  S^{-1}\right)  _{ij}+\sum_{k=1;k\neq i}^{N}S_{ik}\left(
S^{-1}\right)  _{kj}%
\]
Introducing the definition (\ref{def_S}) of $S$ yields%
\[
1_{\left\{  i=j\right\}  }=\left(  S^{-1}\right)  _{ij}\frac{\delta_{i}%
}{\left(  1-v_{i\infty}\right)  ^{2}}-\sum_{k=1;k\neq i}^{N}a_{ik}\beta
_{k}\left(  S^{-1}\right)  _{kj}%
\]
where now the infection rates $\beta_{k}$ need to stay inside the summation.
Rewritten yields the second form for%
\begin{equation}
\left(  S^{-1}\right)  _{ij}=\frac{\left(  1-v_{i\infty}\right)  ^{2}}%
{\delta_{i}}\left(  1_{\left\{  i=j\right\}  }+\sum_{k=1;k\neq i}^{N}%
a_{ik}\beta_{k}\left(  S^{-1}\right)  _{kj}\right)  \label{inverseS_ij_vorm2}%
\end{equation}
},
\begin{align*}
1_{\left\{  i=j\right\}  }  &  =\sum_{k=1}^{N}\left(  S^{-1}\right)
_{ik}S_{kj}\\
&  =\left(  S^{-1}\right)  _{ij}S_{jj}+\sum_{k=1;k\neq j}^{N}\left(
S^{-1}\right)  _{ik}S_{kj}%
\end{align*}
Introducing the definition (\ref{def_S}) of $S$ yields%
\[
1_{\left\{  i=j\right\}  }=\left(  S^{-1}\right)  _{ij}\frac{\delta_{j}%
}{\left(  1-v_{j\infty}\right)  ^{2}}-\sum_{k=1;k\neq j}^{N}\left(
S^{-1}\right)  _{ik}a_{kj}\beta_{j}%
\]
Thus, if $j=i$, then%
\begin{equation}
1=\frac{\delta_{i}\left(  S^{-1}\right)  _{ii}}{\left(  1-v_{i\infty}\right)
^{2}}-\beta_{i}\sum_{k=1;k\neq i}^{N}\left(  S^{-1}\right)  _{ik}a_{ki}
\label{1_from_inverse}%
\end{equation}
from which
\begin{equation}
1\leq\frac{\delta_{i}\left(  S^{-1}\right)  _{ii}}{\left(  1-v_{i\infty
}\right)  ^{2}} \label{inequality_inverse}%
\end{equation}
follows.

We can also write%
\begin{equation}
\left(  S^{-1}\right)  _{ij}=\frac{\left(  1-v_{j\infty}\right)  ^{2}}%
{\delta_{j}}\left(  1_{\left\{  i=j\right\}  }+\beta_{j}\sum_{k=1;k\neq j}%
^{N}\left(  S^{-1}\right)  _{ik}a_{kj}\right)  \label{inverseS_ij_vorm1}%
\end{equation}
from which we find that%
\[
\left(  S^{-1}\right)  _{ij}\leq\frac{\left(  1-v_{j\infty}\right)  ^{2}%
}{\delta_{j}}\left(  1_{\left\{  i=j\right\}  }+\beta_{j}d_{j}\max_{1\leq
k\leq N}\left(  S^{-1}\right)  _{ik}\right)
\]
as well as the lower bound%
\[
\left(  S^{-1}\right)  _{ij}\geq\frac{\left(  1-v_{j\infty}\right)  ^{2}%
}{\delta_{j}}\left(  1_{\left\{  i=j\right\}  }+\beta_{j}d_{j}\min_{1\leq
k\leq N}\left(  S^{-1}\right)  _{ik}\right)
\]

\onecolumn The H\"{o}lder inequality \cite[p. 107]{PVM_PAComplexNetsCUP} with
$p>1$ and $\frac{1}{p}+\frac{1}{q}=1$ shows, using $a_{kj}=a_{kj}^{2}$, that%
\[
\sum_{k=1;k\neq j}^{N}\left(  S^{-1}\right)  _{ik}a_{kj}=\sum_{k=1;k\neq
j}^{N}\left\{  \left(  S^{-1}\right)  _{ik}a_{kj}\right\}  a_{kj}\leq\left(
\sum_{k=1;k\neq j}^{N}a_{kj}\left(  S^{-1}\right)  _{ik}^{p}\right)
^{\frac{1}{p}}\left(  \sum_{k=1;k\neq j}^{N}a_{kj}^{q}\right)  ^{\frac{1}{q}%
}=d_{j}^{1-\frac{1}{p}}\left(  \sum_{k=1;k\neq j}^{N}a_{kj}\left(
S^{-1}\right)  _{ik}^{p}\right)  ^{\frac{1}{p}}%
\]
so that, for $p\geq1$ (including $p=1$ for which equality holds)
\[
\left(  S^{-1}\right)  _{ij}\leq\frac{\left(  1-v_{j\infty}\right)  ^{2}%
}{\delta_{j}}\left(  1_{\left\{  i=j\right\}  }+\beta_{j}d_{j}^{1-\frac{1}{p}%
}\left(  \sum_{k=1;k\neq j}^{N}a_{kj}\left(  S^{-1}\right)  _{ik}^{p}\right)
^{\frac{1}{p}}\right)
\]

For $i=j$ in (\ref{inverseS_ij_vorm1}), we have%
\begin{equation}
\left(  S^{-1}\right)  _{ii}=\left(  1-v_{i\infty}\right)  ^{2}\tau_{i}\left(
\frac{1}{\beta_{i}}+\sum_{k=1;k\neq i}^{N}\left(  S^{-1}\right)  _{ik}%
a_{ki}\right)  \label{identity_from_inverse_matrix}%
\end{equation}
while, if $j\neq i$, then%
\begin{align*}
\left(  S^{-1}\right)  _{ij}  &  =\left(  1-v_{j\infty}\right)  ^{2}\tau
_{j}\sum_{k=1;k\neq j}^{N}\left(  S^{-1}\right)  _{ik}a_{kj}\\
&  =\left(  1-v_{j\infty}\right)  ^{2}\tau_{j}\left(  S^{-1}\right)
_{ii}a_{ij}+\left(  1-v_{j\infty}\right)  ^{2}\tau_{j}\sum_{k=1;k\neq\left\{
i,j\right\}  }^{N}\left(  S^{-1}\right)  _{ik}a_{kj}%
\end{align*}
which illustrates that $\left(  S^{-1}\right)  _{ii}=O\left\{  \frac{\left(
1-v_{i\infty}\right)  ^{2}}{\delta_{i}}\right\}  $ and $\left(  S^{-1}\right)
_{ij}=O\left(  \left(  1-v_{j\infty}\right)  ^{2}\tau_{j}\frac{\left(
1-v_{i\infty}\right)  ^{2}}{\delta_{i}}\right)  $ as $v_{i\infty}\rightarrow
1$. Also, that for $j\neq i$,%
\begin{equation}
\left(  S^{-1}\right)  _{ij}\geq\left(  1-v_{j\infty}\right)  ^{2}\tau
_{j}\left(  S^{-1}\right)  _{ii}a_{ij} \label{lower_bound_S_ij}%
\end{equation}
This inequality can be slightly generalized. Indeed, from
(\ref{inverseS_ij_vorm1}), we have%
\[
\left(  S^{-1}\right)  _{ik}a_{kj}\leq\sum_{k=1;k\neq j}^{N}\left(
S^{-1}\right)  _{ik}a_{kj}=\frac{\left(  S^{-1}\right)  _{ij}}{\tau_{j}\left(
1-v_{j\infty}\right)  ^{2}}-\frac{1_{\left\{  i=j\right\}  }}{\beta_{j}}%
\]
so that, for $k=i$,%
\[
\left(  S^{-1}\right)  _{ij}\geq\tau_{j}\left(  1-v_{j\infty}\right)
^{2}\left(  S^{-1}\right)  _{ii}a_{ij}+\frac{\left(  1-v_{j\infty}\right)
^{2}}{\delta_{j}}1_{\left\{  i=j\right\}  }%
\]
If $j=i$, then $\left(  S^{-1}\right)  _{ii}\geq\frac{\left(  1-v_{i\infty
}\right)  ^{2}}{\delta_{i}}$ and we find (\ref{inequality_inverse}) again.
Similarly, from (\ref{inverseS_ij_vorm2}), we find that%
\[
\left(  S^{-1}\right)  _{ij}\geq\frac{\beta_{j}\left(  1-v_{i\infty}\right)
^{2}}{\delta_{i}}a_{ij}\left(  S^{-1}\right)  _{jj}+\frac{\left(
1-v_{i\infty}\right)  ^{2}}{\delta_{i}}1_{\left\{  i=j\right\}  }%
\]
Thus,%
\begin{equation}
\left(  S^{-1}\right)  _{ij}\geq a_{ij}\max\left(  \left(  1-v_{j\infty
}\right)  ^{2}\tau_{j}\left(  S^{-1}\right)  _{ii},\frac{\beta_{j}}{\delta
_{i}}\left(  1-v_{i\infty}\right)  ^{2}\left(  S^{-1}\right)  _{jj}\right)
\label{lower_bound_S_ij_generalized}%
\end{equation}

Finally, combining the inequality (\ref{inequality_steady_state}) and
(\ref{inequality_inverse}) yields the bounds%
\[
\frac{\left(  1-v_{i\infty}\right)  ^{2}}{\delta_{i}}\leq\left(
S^{-1}\right)  _{ii}<\frac{1}{v_{i\infty}}\frac{\left(  1-v_{i\infty}\right)
^{2}}{\delta_{i}}%
\]
Since $S$ and $S^{-1}$ are positive definite, it holds \cite[p. 241]%
{PVM_graphspectra} that%
\begin{equation}
\left(  S^{-1}\right)  _{ij}\leq\min\left(  \frac{\left(  S^{-1}\right)
_{ii}+\left(  S^{-1}\right)  _{jj}}{2},\sqrt{\left(  S^{-1}\right)
_{ii}\left(  S^{-1}\right)  _{jj}}\right)  \label{upperbound_inverseS_ij}%
\end{equation}

\subsection{Deductions from the the steady-state equation
(\ref{steady_state_equation})}

\label{sec_S_steady_state}We rewrite the steady-state equation
(\ref{steady_state_equation})
\begin{align*}
A\text{diag}\left(  \beta_{i}\right)  V_{\infty}  &  =\text{diag}\left(
\frac{\delta_{i}}{1-v_{i\infty}}\right)  V_{\infty}=\text{diag}\left(
\frac{\delta_{i}}{\left(  1-v_{i\infty}\right)  ^{2}}\right)  \text{diag}%
\left(  1-v_{i\infty}\right)  V_{\infty}\\
&  =\text{diag}\left(  \frac{\delta_{i}}{\left(  1-v_{i\infty}\right)  ^{2}%
}\right)  V_{\infty}-\text{diag}\left(  \frac{\delta_{i}v_{i\infty}}{\left(
1-v_{i\infty}\right)  ^{2}}\right)  V_{\infty}%
\end{align*}
in terms of the matrix $S$ in (\ref{def_S}),%
\[
SV_{\infty}=\text{diag}\left(  \frac{\delta_{i}v_{i\infty}}{\left(
1-v_{i\infty}\right)  ^{2}}\right)  V_{\infty}%
\]
Since the inverse $S^{-1}$ exists above the critical threshold, we arrive at%
\begin{equation}
V_{\infty}=S^{-1}\text{diag}\left(  \frac{\delta_{i}v_{i\infty}}{\left(
1-v_{i\infty}\right)  ^{2}}\right)  V_{\infty}
\label{steady_state_in_terms_of_S}%
\end{equation}
The $i$-th row is%
\begin{align}
v_{i\infty}  &  =\sum_{j=1}^{N}\frac{\delta_{j}v_{j\infty}^{2}\left(
S^{-1}\right)  _{ij}}{\left(  1-v_{j\infty}\right)  ^{2}}%
\label{steady_state_equation_in_terms_of_S_row_i_vi}\\
&  =\frac{\delta_{i}v_{i\infty}^{2}\left(  S^{-1}\right)  _{ii}}{\left(
1-v_{i\infty}\right)  ^{2}}+\sum_{k=1;k\neq i}^{N}\left(  S^{-1}\right)
_{ik}\frac{\delta_{k}v_{k\infty}^{2}}{\left(  1-v_{k\infty}\right)  ^{2}%
}\nonumber
\end{align}
Thus,%
\begin{equation}
1=\frac{\delta_{i}v_{i\infty}\left(  S^{-1}\right)  _{ii}}{\left(
1-v_{i\infty}\right)  ^{2}}+\sum_{k=1;k\neq i}^{N}\left(  S^{-1}\right)
_{ik}\frac{\delta_{k}v_{k\infty}^{2}}{v_{i\infty}\left(  1-v_{k\infty}\right)
^{2}} \label{steady_state_equation_in_terms_of_S_row_i}%
\end{equation}

\subsection{Analysis of $M_{ki}$ defined in (\ref{def_M_ki})}

\label{sec_analysis_M_ki}The results presented in this section illustrate the
difficulty to determine the sign of $M_{ki}$, which prevents us to draw
conclusions about convexity or concavity of $v_{k\infty}$ as a function of
$\delta_{i}$, given that all $\delta_{k}$ are independent.

\textbf{A.} We can write%
\[
\sum_{j=1}^{N}\left(  S^{-1}\right)  _{kj}\frac{\delta_{j}\left(  \left(
S^{-1}\right)  _{ji}\right)  ^{2}}{\left(  1-v_{j\infty}\right)  ^{3}}%
=\sum_{j=1}^{N}\left(  S^{-1}\right)  _{kj}\frac{\delta_{j}\left(
S^{-1}\right)  _{ji}}{\left(  1-v_{j\infty}\right)  ^{3}}\left(
S^{-1}\right)  _{ji}=\left(  S^{-1}\text{diag}\left(  \frac{\delta_{j}\left(
S^{-1}\right)  _{ji}}{\left(  1-v_{j\infty}\right)  ^{3}}\right)
S^{-1}\right)  _{ki}%
\]
Iterating (\ref{steady_state_in_terms_of_S}) once yields%
\begin{equation}
V_{\infty}=S^{-1}\text{diag}\left(  \frac{\delta_{i}v_{i\infty}}{\left(
1-v_{i\infty}\right)  ^{2}}\right)  S^{-1}\text{diag}\left(  \frac{\delta
_{i}v_{i\infty}}{\left(  1-v_{i\infty}\right)  ^{2}}\right)  V_{\infty}
\label{steady_states_in_terms_of_S_iterated_twice}%
\end{equation}
and%
\[
\left(  S^{-1}\text{diag}\left(  \frac{\delta_{i}v_{i\infty}}{\left(
1-v_{i\infty}\right)  ^{2}}\right)  S^{-1}\right)  _{ij}=\sum_{l=1}^{N}\left(
S^{-1}\right)  _{il}\frac{\delta_{l}v_{l\infty}}{\left(  1-v_{l\infty}\right)
^{2}}\left(  S^{-1}\right)  _{lj}%
\]
The corresponding $i$-th row in
(\ref{steady_states_in_terms_of_S_iterated_twice}) is%

\begin{equation}
v_{i\infty}=\sum_{j=1}^{N}\sum_{l=1}^{N}\left(  S^{-1}\right)  _{il}%
\frac{\delta_{l}v_{l\infty}}{\left(  1-v_{l\infty}\right)  ^{2}}\left(
S^{-1}\right)  _{lj}\frac{\delta_{j}v_{j\infty}^{2}}{\left(  1-v_{j\infty
}\right)  ^{2}} \label{vi_infty_steady_state_twice_iterated}%
\end{equation}
which illustrates that the double sum containing products of elements of
$S^{-1}$ can be smaller than 1. In addition, using
(\ref{steady_state_equation_in_terms_of_S_row_i_vi}) into
(\ref{vi_infty_steady_state_twice_iterated}) agains leads to
(\ref{steady_state_equation_in_terms_of_S_row_i_vi}).

\textbf{B.} Since all elements of $S^{-1}$ are positive, we have that%
\[
\left(  S^{-2}\right)  _{ki}=\sum_{j=1}^{N}\left(  S^{-1}\right)  _{kj}\left(
S^{-1}\right)  _{ji}\geq\left(  S^{-1}\right)  _{ki}\left(  S^{-1}\right)
_{ii}%
\]
which we use to lower bound $M_{ki}$ as%
\begin{align*}
M_{ki}  &  =\frac{\left(  S^{-1}\right)  _{ki}\left(  S^{-1}\right)  _{ii}%
}{\left(  1-v_{i\infty}\right)  }-v_{i\infty}\sum_{j=1}^{N}\left(
S^{-1}\right)  _{kj}\frac{\delta_{j}\left(  \left(  S^{-1}\right)
_{ji}\right)  ^{2}}{\left(  1-v_{j\infty}\right)  ^{3}}\\
&  \leq\frac{1}{\left(  1-v_{i\infty}\right)  }\sum_{j=1}^{N}\left(
S^{-1}\right)  _{kj}\left(  S^{-1}\right)  _{ji}-v_{i\infty}\sum_{j=1}%
^{N}\left(  S^{-1}\right)  _{kj}\frac{\delta_{j}\left(  S^{-1}\right)  _{ji}%
}{\left(  1-v_{j\infty}\right)  ^{3}}\left(  S^{-1}\right)  _{ji}\\
&  =\sum_{j=1}^{N}\left(  S^{-1}\right)  _{kj}\left\{  \frac{1}{\left(
1-v_{i\infty}\right)  }-\frac{v_{i\infty}\delta_{j}\left(  S^{-1}\right)
_{ji}}{\left(  1-v_{j\infty}\right)  ^{3}}\right\}  \left(  S^{-1}\right)
_{ji}%
\end{align*}
The equation (\ref{steady_state_equation_in_terms_of_S_row_i_vi}), rewritten
as $v_{j\infty}=\sum_{k=1}^{N}\frac{\delta_{k}v_{k\infty}^{2}\left(
S^{-1}\right)  _{jk}}{\left(  1-v_{k\infty}\right)  ^{2}}$, shows that%
\[
v_{j\infty}\geq\frac{\delta_{i}v_{i\infty}^{2}\left(  S^{-1}\right)  _{ji}%
}{\left(  1-v_{i\infty}\right)  ^{2}}\text{ or }v_{j\infty}\frac{\left(
1-v_{i\infty}\right)  ^{2}}{\delta_{i}v_{i\infty}^{2}}\geq\left(
S^{-1}\right)  _{ji}%
\]
so that%
\[
0\leq\frac{1}{1-v_{i\infty}}-f_{ij}\frac{v_{i\infty}\delta_{j}\left(
S^{-1}\right)  _{ji}}{\left(  1-v_{j\infty}\right)  ^{3}}%
\]
with%
\[
f_{ij}=\frac{\delta_{i}v_{i\infty}\left(  1-v_{j\infty}\right)  ^{3}}%
{\delta_{j}v_{j\infty}\left(  1-v_{i\infty}\right)  ^{3}}%
\]
The terms in the sum in the above inequality for $M_{ki}$ is positive if
$f_{ij}\leq1$. Since we cannot show that for all $j$, it holds that $f_{ij}$,
we cannot conclude that the upper bound is always positive.

\textbf{C.} Starting from (\ref{steady_state_equation_in_terms_of_S_row_i_vi})
and assuming that $S$ is symmetric (which happens if all infection rates
$\beta_{k}=\beta$ are the same) so that $\left(  S^{-1}\right)  _{ji}=\left(
S^{-1}\right)  _{ij}$, the Cauchy-Schwarz inequality \cite[p. 107]%
{PVM_PAComplexNetsCUP} shows that%

\begin{align*}
v_{i\infty}^{2}  &  =\left(  \sum_{j=1}^{N}\frac{\delta_{j}v_{j\infty}%
^{2}\left(  S^{-1}\right)  _{ij}}{\left(  1-v_{j\infty}\right)  ^{2}}\right)
^{2}=\left(  \sum_{j=1}^{N}\frac{v_{j\infty}^{2}\sqrt{\delta_{j}}}%
{\sqrt{\left(  S^{-1}\right)  _{kj}\left(  1-v_{j\infty}\right)  }}\frac
{\sqrt{\left(  S^{-1}\right)  _{kj}\delta_{j}}\left(  S^{-1}\right)  _{ij}%
}{\left(  1-v_{j\infty}\right)  ^{3/2}}\right)  ^{2}\\
&  \leq\sum_{j=1}^{N}\frac{v_{j\infty}^{4}\delta_{j}}{\left(  1-v_{j\infty
}\right)  \left(  S^{-1}\right)  _{kj}}\sum_{j=1}^{N}\frac{\left(
S^{-1}\right)  _{kj}\delta_{j}\left(  \left(  S^{-1}\right)  _{ij}\right)
^{2}}{\left(  1-v_{j\infty}\right)  ^{3}}%
\end{align*}
Hence,%
\begin{align*}
v_{i\infty}\sum_{j=1}^{N}\left(  S^{-1}\right)  _{kj}\frac{\delta_{j}\left(
\left(  S^{-1}\right)  _{ij}\right)  ^{2}}{\left(  1-v_{j\infty}\right)
^{3}}  &  \geq\frac{v_{i\infty}^{3}}{\sum_{j=1}^{N}\frac{v_{j\infty}^{4}%
\delta_{j}}{\left(  1-v_{j\infty}\right)  \left(  S^{-1}\right)  _{kj}}}%
=\frac{v_{i\infty}^{3}}{\frac{v_{i\infty}^{4}\delta_{i}}{\left(  1-v_{i\infty
}\right)  \left(  S^{-1}\right)  _{ki}}+\sum_{j=1;j\neq i}^{N}\frac
{v_{j\infty}^{4}\delta_{j}}{\left(  1-v_{j\infty}\right)  \left(
S^{-1}\right)  _{kj}}}\\
&  =\frac{\left(  1-v_{i\infty}\right)  \left(  S^{-1}\right)  _{ki}%
}{v_{i\infty}\delta_{i}\left(  \sum_{j=1}^{N}\frac{v_{j\infty}^{4}\delta
_{j}\left(  1-v_{i\infty}\right)  \left(  S^{-1}\right)  _{ki}}{v_{i\infty
}^{4}\delta_{i}\left(  1-v_{j\infty}\right)  \left(  S^{-1}\right)  _{kj}%
}\right)  }%
\end{align*}
so that%
\begin{align*}
\left(  -M_{ki}\right)   &  \geq\frac{\left(  1-v_{i\infty}\right)  \left(
S^{-1}\right)  _{ki}}{v_{i\infty}\delta_{i}\left(  \sum_{j=1}^{N}%
\frac{v_{j\infty}^{4}\delta_{j}\left(  1-v_{i\infty}\right)  \left(
S^{-1}\right)  _{ki}}{v_{i\infty}^{4}\delta_{i}\left(  1-v_{j\infty}\right)
\left(  S^{-1}\right)  _{kj}}\right)  }-\frac{\left(  S^{-1}\right)
_{ii}\left(  S^{-1}\right)  _{ki}}{\left(  1-v_{i\infty}\right)  }\\
&  =\left(  S^{-1}\right)  _{ki}\left\{  \frac{\left(  1-v_{i\infty}\right)
}{v_{i\infty}\delta_{i}\left(  \sum_{j=1}^{N}\frac{v_{j\infty}^{4}\delta
_{j}\left(  1-v_{i\infty}\right)  \left(  S^{-1}\right)  _{ki}}{v_{i\infty
}^{4}\delta_{i}\left(  1-v_{j\infty}\right)  \left(  S^{-1}\right)  _{kj}%
}\right)  }-\frac{\left(  S^{-1}\right)  _{ii}}{\left(  1-v_{i\infty}\right)
}\right\}
\end{align*}
Invoking (\ref{identity_from_inverse_matrix}) yields%
\begin{align*}
\left(  -M_{ki}\right)   &  \geq\left(  S^{-1}\right)  _{ki}\left\{
\frac{\left(  1-v_{i\infty}\right)  }{v_{i\infty}\delta_{i}\left(  \sum
_{j=1}^{N}\frac{v_{j\infty}^{4}\delta_{j}\left(  1-v_{i\infty}\right)  \left(
S^{-1}\right)  _{ki}}{v_{i\infty}^{4}\delta_{i}\left(  1-v_{j\infty}\right)
\left(  S^{-1}\right)  _{kj}}\right)  }-\frac{\left(  1-v_{i\infty}\right)
}{\delta_{i}}-\left(  1-v_{i\infty}\right)  \tau_{i}\sum_{l=1;l\neq i}%
^{N}\left(  S^{-1}\right)  _{il}a_{li}\right\} \\
&  =\left(  S^{-1}\right)  _{ki}\frac{\left(  1-v_{i\infty}\right)  }%
{\delta_{i}}\left\{  \frac{1}{v_{i\infty}\left(  1+\sum_{j=1;j\neq i}^{N}%
\frac{v_{j\infty}^{4}\delta_{j}\left(  1-v_{i\infty}\right)  \left(
S^{-1}\right)  _{ki}}{v_{i\infty}^{4}\delta_{i}\left(  1-v_{j\infty}\right)
\left(  S^{-1}\right)  _{kj}}\right)  }-\left(  1+\beta_{i}\sum_{l=1;l\neq
i}^{N}\left(  S^{-1}\right)  _{il}a_{li}\right)  \right\}
\end{align*}
Combining the key inequality (\ref{inequality_steady_state}), which is
equivalent to%
\[
\frac{1}{v_{i\infty}}>\frac{\delta_{i}\left(  S^{-1}\right)  _{ii}}{\left(
1-v_{i\infty}\right)  ^{2}}%
\]
with (\ref{identity_from_inverse_matrix}), leads to%
\[
\frac{1}{v_{i\infty}\left(  1+\beta_{i}\sum_{l=1;l\neq i}^{N}\left(
S^{-1}\right)  _{il}a_{li}\right)  }>1
\]
The condition for $\left(  -M_{ki}\right)  $ to be positive is%
\[
\frac{1}{v_{i\infty}\left(  1+\sum_{j=1;j\neq i}^{N}\frac{v_{j\infty}%
^{4}\delta_{j}\left(  1-v_{i\infty}\right)  \left(  S^{-1}\right)  _{ki}%
}{v_{i\infty}^{4}\delta_{i}\left(  1-v_{j\infty}\right)  \left(
S^{-1}\right)  _{kj}}\right)  }\geq\left(  1+\beta_{i}\sum_{l=1;l\neq i}%
^{N}\left(  S^{-1}\right)  _{il}a_{li}\right)
\]
or%
\[
\frac{1}{v_{i\infty}\left(  1+\beta_{i}\sum_{l=1;l\neq i}^{N}\left(
S^{-1}\right)  _{il}a_{li}\right)  }\geq1+\frac{\left(  1-v_{i\infty}\right)
\left(  S^{-1}\right)  _{ki}}{v_{i\infty}^{4}\delta_{i}}\sum_{j=1;j\neq i}%
^{N}\frac{v_{j\infty}^{4}\delta_{j}}{\left(  1-v_{j\infty}\right)  \left(
S^{-1}\right)  _{kj}}%
\]
Again, in general, this condition is difficult to assess and there might be a
region for $\delta_{i}$ (or $v_{k\infty}$) where the condition is satisfied
(thus, where $v_{k\infty}$ is concave in $\delta_{i}$).

\textbf{D.} When introducing (\ref{inverseS_ij_vorm2})%
\[
\left(  S^{-1}\right)  _{ji}=\frac{\left(  1-v_{j\infty}\right)  ^{2}}%
{\delta_{j}}\left(  1_{\left\{  i=j\right\}  }+\sum_{k=1;k\neq j}^{N}%
a_{jk}\beta_{k}\left(  S^{-1}\right)  _{ki}\right)
\]
into\footnote{Substitution of (\ref{inverseS_ij_vorm1}) leads to less
transparent equations.} $M_{ki}$, we obtain%
\begin{align*}
M_{ki} &  =\frac{\left(  S^{-1}\right)  _{ki}\left(  S^{-1}\right)  _{ii}%
}{\left(  1-v_{i\infty}\right)  }-v_{i\infty}\sum_{j=1}^{N}\left(
S^{-1}\right)  _{kj}\frac{\delta_{j}\left(  \left(  S^{-1}\right)
_{ji}\right)  ^{2}}{\left(  1-v_{j\infty}\right)  ^{3}}\\
&  =\frac{\left(  S^{-1}\right)  _{ki}\left(  S^{-1}\right)  _{ii}}{\left(
1-v_{i\infty}\right)  }-v_{i\infty}\sum_{j=1}^{N}\left(  S^{-1}\right)
_{kj}\frac{\left(  S^{-1}\right)  _{ji}}{\left(  1-v_{j\infty}\right)
}\left(  1_{\left\{  i=j\right\}  }+\sum_{l=1;l\neq j}^{N}a_{jl}\beta
_{l}\left(  S^{-1}\right)  _{li}\right)  \\
&  =\frac{\left(  S^{-1}\right)  _{ki}\left(  S^{-1}\right)  _{ii}}{\left(
1-v_{i\infty}\right)  }-v_{i\infty}\left(  S^{-1}\right)  _{ki}\frac{\left(
S^{-1}\right)  _{ii}}{\left(  1-v_{i\infty}\right)  }-v_{i\infty}\sum
_{j=1}^{N}\left(  S^{-1}\right)  _{kj}\frac{\left(  S^{-1}\right)  _{ji}%
}{\left(  1-v_{j\infty}\right)  }\sum_{l=1;l\neq j}^{N}a_{jl}\beta_{l}\left(
S^{-1}\right)  _{li}%
\end{align*}
Thus,%
\begin{align*}
M_{ki} &  =\left(  S^{-1}\right)  _{ki}\left(  S^{-1}\right)  _{ii}%
-v_{i\infty}\sum_{j=1}^{N}\left(  S^{-1}\right)  _{kj}\frac{\left(
S^{-1}\right)  _{ji}}{\left(  1-v_{j\infty}\right)  }\sum_{l=1;l\neq j}%
^{N}a_{jl}\beta_{l}\left(  S^{-1}\right)  _{li}\\
&  =\left(  S^{-1}\right)  _{ki}\left(  S^{-1}\right)  _{ii}\left\{
1-\frac{v_{i\infty}}{\left(  1-v_{i\infty}\right)  }\sum_{l=1;l\neq i}%
^{N}a_{il}\beta_{l}\left(  S^{-1}\right)  _{li}\right\}  -v_{i\infty}%
\sum_{j=1;j\neq i}^{N}\left(  S^{-1}\right)  _{kj}\frac{\left(  S^{-1}\right)
_{ji}}{\left(  1-v_{j\infty}\right)  }\sum_{l=1;l\neq j}^{N}a_{jl}\beta
_{l}\left(  S^{-1}\right)  _{li}%
\end{align*}
We concentrate on the term%
\[
G=1-\frac{v_{i\infty}}{\left(  1-v_{i\infty}\right)  }\sum_{l=1;l\neq i}%
^{N}a_{il}\beta_{l}\left(  S^{-1}\right)  _{li}<1
\]
and substitute (\ref{inverseS_ij_vorm2})%
\[
\frac{\delta_{i}\left(  S^{-1}\right)  _{ii}}{\left(  1-v_{i\infty}\right)
^{2}}-1=\sum_{l=1;l\neq i}^{N}a_{il}\beta_{l}\left(  S^{-1}\right)  _{li}%
\]
so that%
\begin{align*}
G &  =1-\frac{v_{i\infty}\delta_{i}\left(  S^{-1}\right)  _{ii}}{\left(
1-v_{i\infty}\right)  ^{3}}+\left(  1-\frac{1}{1-v_{i\infty}}\right)  \\
&  =2-\frac{1}{1-v_{i\infty}}\left(  \frac{v_{i\infty}\delta_{i}\left(
S^{-1}\right)  _{ii}}{\left(  1-v_{i\infty}\right)  ^{2}}+1\right)  >2\left(
1-\frac{1}{1-v_{i\infty}}\right)  =-2\frac{v_{i\infty}}{1-v_{i\infty}}%
\end{align*}
Hence,%
\[
-2\frac{v_{i\infty}}{1-v_{i\infty}}<G<1
\]
which indicates that, for small $v_{i\infty}$, $G$ can be negative and, in
absolute value larger than the remaining sum in $M_{ki}$.

\end{document}